\magnification=\magstep1
\input amstex
\UseAMSsymbols
\input pictex
\vsize=23truecm
\NoBlackBoxes
\parindent=18pt

   \font\rmk=cmr8

\def\G{\text{\rm G}}
\def\TR{\text{\rm TR}}
\def\op{\text{\rm op}}

\def\gp{\operatorname{gp}}
\def\mod{\operatorname{mod}}

\def\Hom{\operatorname{Hom}}
\def\End{\operatorname{End}}
\def\Ext{\operatorname{Ext}}

\def\rad{\operatorname{rad}}
\def\add{\operatorname{add}}
\def\Ker{\operatorname{Ker}}
\def\Cok{\operatorname{Cok}}
\def\soc{\operatorname{soc}}
\def\Tr{\operatorname{Tr}}
\def\Im{\operatorname{Im}}

\def\K{\operatorname{K}}

\def\injdim{\operatorname{inj.dim.}}
  \def\ss{\ssize }
\def\arr#1#2{\arrow <1.5mm> [0.25,0.75] from #1 to #2}

\def\s{\hfill \square}  	

\vglue1.5truecm
\centerline{\bf  Gorenstein-projective and}
                     \smallskip
\centerline{\bf  semi-Gorenstein-projective modules}
                     \bigskip
\centerline{Claus Michael Ringel, Pu Zhang}
                \bigskip\medskip

\noindent {\narrower Abstract:  \rmk Let $\ss A$ be an artin algebra.
An $\ss A$-module $\ss M$ will be said to be
semi-Gorenstein-projective provided that
$\ss\Ext^i(M,A) = 0$ for all $\ss i\ge 1.$ All Gorenstein-projective
modules are semi-Gorenstein-projective and only few and quite complicated examples of
semi-Gorenstein-projective modules
which are not Gorenstein-projective have been known.
One of the aims of the paper is to provide conditions on $\ss A$ such
that all semi-Gorenstein-projective left modules are Gorenstein-projective (we call such an algebra
left weakly  Gorenstein). In particular, we show that in case there are only finitely many isomorphism
classes of indecomposable left modules which are both semi-Gorenstein-projective and torsionless, then
$\ss A$ is left weakly  Gorenstein. On the other hand, we exhibit a 6-dimensional algebra $\ss\Lambda$
with a semi-Gorenstein-projective module $\ss M$ which is not torsionless (thus not
Gorenstein-projective). Actually, also the $\ss \Lambda$-dual module $\ss M^*$ is
semi-Gorenstein-projective.  In this way, we show
the independence of the total reflexivity conditions of Avramov and Martsinkovsky,
thus completing a partial proof by Jorgensen and \c Sega. Since all the syzygy-modules of $\ss M$
and $\ss M^*$ are
3-dimensional, the example can be checked (and visualized) quite easily.
\par}
	\medskip
{\rmk Key words and phrases. Gorenstein-projective module, semi-Gorenstein-projective module,
left weakly Gorenstein algebra, torsionless module, reflexive module, $\ss t$-torsionfree module,
Frobenius category, $\ss \mho$-quiver.
	\medskip
Supported by NSFC 11431010
\par}

	\bigskip\medskip
{\bf 1. Introduction.}
	\medskip
{\bf 1.1. Notations and definitions.}
Let $A$ be an artin algebra. All modules will be finitely generated.
Usually, the modules we are
starting with will be left modules, but some constructions then yield right modules.
Let $\mod A$ be the category of all finitely generated left $A$-modules
and $\add(A)$ the full subcategory of all projective modules.

If $M$ is a module, let $PM$ be a projective cover of $M$, and
$\Omega M$ the kernel of the canonical map $PM \to M$. The modules $\Omega^tM$ with $t\ge 0$
are called the {\it syzygy} modules of $M$.
A module $M$ is said to be {\it $\Omega$-periodic} provided that
there is some $t\ge 1$ with $\Omega^tM = M.$

The right $A$-module $M^* = \Hom(M,A)$ is called the {\it $A$-dual} of $M$.
Let $\phi_M\:M \to M^{**}$ be defined by $\phi_M(m)(f) = f(m)$ for $m\in M,\ f\in M^*$.
A module $M$ is said to be {\it torsionless} provided that $M$ is a submodule of a projective
module, or, equivalently, provided that $\phi_M$ is injective.
A module $M$ is called {\it reflexive} provided that $\phi_M$ is bijective.

Let $\Tr M$ be the cokernel of $f^*,$ where $f$ is a minimal projective presentation of $M$
(this is the canonical map $P(\Omega M) \to P M$).
Note that $\Tr M$ is a right $A$-module, called the {\it transpose
of $M$}.

{\it A complete projective resolution} is a (double infinite) exact sequence
$$ P^\bullet: \qquad  \cdots \longrightarrow P^{-1}\longrightarrow P^{0}
  \overset {d^0} \to {\longrightarrow} P^{1}\longrightarrow \cdots$$
of projective left $A$-modules, such that $\Hom_A(P^\bullet,A)$ is again exact.
 A module $M$
is {\it Gorenstein-projective}, if there is a complete
projective resolution $P^\bullet$ with $M$ isomorphic to the image of $d^0.$

A module $M$ will be said to be {\it semi-Gorenstein-projective} provided that $\Ext^i(M,A) = 0$
for all $i\ge 1.$ All Gorenstein-projective modules are semi-Gorenstein-projective.
If $M$ is semi-Gorenstein-projective, then so is $\Omega M$.
Denote by $\gp(A)$ the class of all Gorenstein-projective modules
and by ${}^\perp A$ the class of all semi-Gorenstein-projective modules. Then
$\gp(A) \subseteq {}^\perp A.$ We propose to call an artin algebra $A$ {\it left weakly  Gorenstein}
provided that ${}^\perp A = \gp(A)$, i.e.,
any semi-Gorenstein-projective module is Gorenstein-projective.
(And $A$ is said to be {\it right weakly Gorenstein} if its
opposite algebra $A^{\op}$ is left weakly Gorenstein.)

The first aim of the paper is to provide a systematic treatment of the relationship between
semi-Gorenstein-projective modules and Gorenstein-projective modules, see theorems 1.2 to 1.4.
Some of these results are (at least partially) known or can be obtained from the
literature, in particular see the paper [B3] by Beligiannis, but we hope that a unified,
elementary and direct presentation may be appreciated.
	\medskip
{\bf 1.2.} First, we have various
characterizations of the left weakly  Gorenstein algebras.
	\smallskip
{\bf Theorem.} {\it Let $A$ be an artin algebra. The following statements are equivalent$:$
\item{\rm (1)} $A$ is left weakly  Gorenstein.
\item{\rm (2)} Any semi-Gorenstein-projective module is torsionless.
\item{\rm (3)} Any semi-Gorenstein-projective module is reflexive.
\item{\rm (4)} For any semi-Gorenstein-projective module $M$, the map $\phi_M$ is surjective.
\item{\rm (5)} For any semi-Gorenstein-projective module $M$, the module $M^*$ is semi-Gorenstein projective.
\item{\rm (6)} Any semi-Gorenstein-projective module $M$ satisfies $\Ext^1(M^*,A_A) = 0$.
\item{\rm (7)} Any semi-Gorenstein-projective module $M$ satisfies $\Ext^1(\Tr M,A_A) = 0.$\par}
	\medskip
The equivalence of (1) and (2) was published by Huang-Huang [HH, Theorem 4.2].
	\bigskip
{\bf 1.3.} The next result concerns artin algebras with finitely many indecomposable
semi-Gorenstein-projective modules or with finitely many indecomposable torsionless modules.
	\smallskip
{\bf Theorem.} {\it If the number of isomorphism classes of indecomposable modules which are both
semi-Gorenstein-projective and torsionless is finite, then $A$ is left weakly  Gorenstein
and any indecomposable non-projective semi-Gorenstein-projective module is $\Omega$-periodic.}
	\smallskip
This combines two different directions of thoughts. First of all, Yoshino [Y] has shown
that for certain commutative rings $R$ (in particular all artinian commutative rings) the
finiteness of the number of isomorphism classes of indecomposable
semi-Gorenstein-projective $R$-modules implies that $R$ is left weakly  Gorenstein. This was generalized
to artin algebras by Beligiannis [B3, Corollary 5.11].
Second, according to Marczinzik [M1], all torsionless-finite artin algebras (these are the artin algebras with only
finitely many isomorphism classes of torsionless indecomposable modules) are left weakly  Gorenstein.
Note that a lot of interesting classes of artin algebras are
torsionless-finite, see 3.6.
	\bigskip
{\bf 1.4.} If $\Cal C$ is an extension-closed full subcategory of $\mod A$, then
the embedding of $\Cal C$ into $\mod A$ provides an exact structure on $\Cal C$,
called its {\it canonical} exact structure
(for the basic properties of exact structures, see for example [K, Appendix A]).
An exact category $\Cal F$ is called a {\it Frobenius category} provided
that it has enough
projective and enough injective objects and that the
projective objects in $\Cal F$ are just the injective objects in $\Cal F$.
We denote by $\Cal P(\Cal F)$ (and by $\Cal I(\Cal F)$) the full subcategory of the
projective (respectively injective) objects in $\Cal F$.

The subcategories $\gp(A)$ and ${}^\perp A$ are extension-closed, and with its canonical exact structure $\gp(A)$ is Frobenius with $\Cal P(\gp(A)) = \add A$ ([B2, Prop. 3.8]).
Thus, if $A$ is left weakly Gorenstein, then  $\Cal F = \gp(A)$ is
an  extension-closed full subcategory of $\mod A$ which is Frobenius with the
canonical exact structure and satisfies
$\Cal P(\Cal F) \subseteq {}^\perp A \subseteq \Cal F$.
The following result shows that these properties characterize left weakly Gorenstein algebras.
	\smallskip
{\bf Theorem.} {\it Let $A$ be an artin algebra and $\Cal F$
an  extension-closed full subcategory of $\mod A$ such that
$\Cal F$ is a Frobenius category with respect to its canonical exact structure.
If $\Cal P(\Cal F) \subseteq {}^\perp A \subseteq \Cal F$, then
$A$ is left weakly Gorenstein and
$\Cal F = \gp(A).$}
	\medskip
A full subcategory $\Cal C$ of $\mod A$ is said to be {\it resolving} provided
that it contains all
the projective modules and is closed under extensions, direct summands and kernels
of surjective maps.
Note that ${}^\perp\! A$ and $\gp(A)$ are resolving subcategories.
	\smallskip
{\bf Corollary 1.} {\it Let $A$ be an artin algebra and
$\Cal F$ a resolving subcategory of $\mod A$ with ${}^\perp\! A \subseteq \Cal F$. Assume that
$\Cal F$ with its canonical exact structure is a Frobenius subcategory. Then
$A$ is left weakly Gorenstein and
$\Cal F = \gp(A).$}
	\medskip
Taking $\Cal F = \ ^\perp A$  in Theorem 1.4 we get
	\smallskip
{\bf Corollary 2.} {\it An artin algebra $A$ is left weakly  Gorenstein if and only if
${}^\perp\! A$ with its canonical exact structure is a Frobenius subcategory.}
	\medskip
We remark that $\gp(A)$ is the largest resolving Frobenius subcategory of $\mod A$ (compare [B1, Prop. 2.13, Theorem 2.11], [B2, p.145], and [B3, p.1989]; also [ZX, Prop. 5.1]).
This implies Theorem 1.4 and the two corollaries (as one of the referees has pointed out).
	\bigskip
{\bf 1.5. The $\mho$-quiver of an artin algebra $A$.}
The main tool used in the paper are the operator $\mho$,
and the $\mho$-quiver of $A$. Here are the definitions.

Recall that a map $f\:M \to M'$ is said to be {\it left minimal} provided that
any map $h\:M' \to M'$ with $hf = f$ is an automorphism ([AR1]).
A left $\add(A)$-approximation
will be called {\it minimal} provided that it is left minimal. We denote by
$\mho M$ the cokernel of a minimal left $\add(A)$-approximation of $M$.
(The symbol $\mho$, pronounced
``agemo'', should be a reminder that $\mho$
has to be considered as a kind of inverse of $\Omega$.)
It turns out that the operator $\mho$ coincides with $\Tr\Omega\Tr$, which has been studied by Auslander and Reiten in [AR2]. See Subsection 4.4, and also 4.7.

Let $\omega\:M \to P$ be a minimal left $\add(A)$-approximation with cokernel map
$\pi\:P \to \mho M$. If $M$ is indecomposable and not projective, then the image of $\omega$
is contained in the radical of $P$, thus $\pi$ is a projective cover.
If $M$ is, in addition, torsionless (so that $\omega$ is injective), then $\mho M$ is indecomposable
and not projective, and $\Omega \mho M \simeq M.$

The {\it $\mho$-quiver} of $A$ has as vertices the isomorphism classes $[X]$ of the
indecomposable non-projective modules $X$ and there is an  arrow
$$
{\beginpicture
\setcoordinatesystem units <1cm,1cm>
\put{$[X]$} at 0 0
\setdashes <1mm>
\arr{1.5 0}{0.5 0}
\put{$[\mho X]$} at 2 0
\endpicture}
$$
for any torsionless (indecomposable, non-projective) module $X$.
(We hope that the reader is not irritated
by the chosen orientation of the arrow:
it corresponds to the usual convention when dealing with $\Ext$-quivers.)
A component of the $\mho$-quiver will be called an {\it $\mho$-component;}
a path in the $\mho$-quiver will be called an {\it $\mho$-path.}
	
In the $\mho$-quiver, an arrow ending at $[X]$ starts at $[\mho X]$, thus for any vertex $[X]$, there
is at most one arrow ending in $[X]$. If $[Z]$ is the start of an arrow, say $Z \simeq \mho X$ for
some vertex $[X]$, then $X\simeq \Omega\mho X \simeq \Omega Z$ implies that the arrow is uniquely
determined. This shows that
{\it at any vertex of the $\mho$-quiver, at most one arrow starts and at most one arrow
ends.} As a consequence, we have:
	\medskip
{\bf Proposition.} {\it Any $\mho$-component is a linearly oriented quiver
$\Bbb A_n$ with $n\ge 1$ vertices, or an oriented cycle $\widetilde{\Bbb A}_n$ with
$n\!+\!1\ge 1$ vertices, or of the form $-\Bbb N,$ or $\Bbb N,$  or $\Bbb Z$.}
	\medskip
Note that we consider any subsets $I$ of $\Bbb Z$
as a quiver, with arrows from $z$ to $z\!-\!1$
(provided that both $z\!-\!1$ and $z$ belong to $I$). For example, the interval
$\{1,2,\dots,n\}$ is the quiver $\Bbb A_n$ with linear orientation (with $1$ being
the unique sink and $n$ the unique source).
Here are the quivers $-\Bbb N$ and $\Bbb N$:
$$
{\beginpicture
\setcoordinatesystem units <1cm,1cm>
\put{\beginpicture
\multiput{$\circ$} at 1 0  2 0  3  0 /
\setdashes <1mm>
\arr{0.8 0}{0.2 0}
\arr{1.8 0}{1.2 0}
\arr{2.8 0}{2.2 0}
\put{$\cdots$} at -.3 0
\put{$-\Bbb N$} at 1.5 -.5
\endpicture} at 0 0
\put{\beginpicture
\multiput{$\circ$} at 0 0  1 0  2 0 /
\setdashes <1mm>
\arr{0.8 0}{0.2 0}
\arr{1.8 0}{1.2 0}
\arr{2.8 0}{2.2 0}
\put{$\cdots$} at 3.3 0
\put{$\Bbb N$} at 1.5 -.5
\endpicture} at 5 0
\endpicture}
$$

As we will see in 7.7, all cases mentioned here can arise as $\mho$-components.
	\smallskip
An indecomposable non-projective module $M$ will be said to be of {\it $\mho$-type}
$\Delta$ where
$\Delta \in \{\Bbb A_n,\ \widetilde{\Bbb A}_n,\, -\Bbb N,\ \Bbb N,\ \Bbb Z\}$
in case the $\mho$-component containing $[M]$ is of the form $\Delta$.
	\smallskip
Let us collect what can be read out about an indecomposable non-projective
module when looking at its position in the $\mho$-quiver.
Recall that a module $M$ is said to be
{\it $t$-torsionfree,} provided $\Ext^i(\Tr M,A_A) = 0$
for $1\le i \le t$ (and {\it $\infty$-torsionfree,} provided $\Ext^i(\Tr M,A_A) = 0$
for all $i \ge 1$); the definition is due to Auslander (see [A1, Br, AB]).
	\medskip
{\bf Theorem.}
{\it Let $M$ be an indecomposable non-projective module.
	\smallskip
\item{\rm(0)} $[M]$ is an isolated vertex \ iff \ $\Ext^1(M,A) \neq 0$ and $M$ is not torsionless.
	\smallskip
\item{\rm(1)} $[M]$ is the start of a path of length $t\ge 1$ \ iff \
  $\Ext^i(M,A) = 0$ for $1\le i \le t.$ \newline
  In particular, $[M]$ is the start of an arrow \ iff \  $\Ext^1(M,A) = 0.$

\item{\rm(1$'$)} $[M]$ is the start of an infinite path \ iff \  $M$ is semi-Gorenstein-projective.
\item{\rm(1$''$)} $[M]$ is of $\mho$-type $-\Bbb N$ \ iff \ $M$ is semi-Gorenstein-projective, but
   not Gorenstein-projective.
 \medskip
\item{\rm(2)} $[M]$ is the end of a path of length $t\ge 1$ \ iff \
  $M$ is $t$-torsionfree  for $1\le i \le t$,  \ iff \  $\mho^{i-1} M$ is torsionless for $1\le i \le t$. \newline
 In particular,  $[M]$ is the end of an arrow \ iff \  $M$ is torsionless$;$
 and  $[M]$ is the end of a path of length $2$ \ iff \  $M$ is reflexive.

\item{\rm(2$'$)}  $[M]$ is the end of an infinite path \ iff \
  $M$ is $\infty$-torsionfree, \ iff \
 $M$ is reflexive and $M^*$ is semi-Gorenstein-projective.
\item{\rm(2$''$)} $[M]$ is of $\mho$-type $\Bbb N$ \ iff \ $M$ is $\infty$-torsionfree, but
   not Gorenstein-projective.
 \medskip
\item{\rm(3)} $[M]$ is the start of an infinite path and also the end of an infinite path
\ iff \  $M$ is Gorenstein-projective. \newline
$M$ is of $\mho$-type $\Bbb Z$ \ iff \  $M$ is Gorenstein-projective
 and not $\Omega$-periodic.
 \newline
$M$  is of $\mho$-type $\widetilde{\Bbb A}_n$ for some $n\ge 0$\ iff \
    $M$ is Gorenstein-projective and $\Omega$-periodic.
	\medskip
\item{\rm(4)} $A$-duality provides a bijection between the isomorphism classes of the
  reflexive indecomposable $A$-modules
  of $\mho$-type $\Bbb A_n$ and the isomorphism classes of the reflexive indecomposable $A^{\op}$-modules
  of $\mho$-type $\Bbb A_n$. \newline
  Thus, for any $n\ge 3$,
$A$ has $\mho$-components of form $\Bbb A_n$ \ iff \
$A^{\op}$ has $\mho$-components of form $\Bbb A_n$.
	\smallskip
\item{\rm(5)}
 $A$-duality provides a bijection between the isomorphism classes of the
  reflexive indecomposable $A$-modules
  of $\mho$-type $\Bbb N$ and the isomorphism classes of the reflexive indecomposable $A^{\op}$-modules
  of $\mho$-type $-\Bbb N$.
  \newline
Thus, $A$ has $\mho$-components of form $\Bbb N$ \ iff \
$A^{\op}$ has $\mho$-components of form $-\Bbb N$.
 \par}

	\medskip
{\bf Remark 1. Characterizations of Gorenstein-projective modules.}
The $\mho$-quiver shows nicely that an indecomposable module $M$ is
Gorenstein-projective if and only if both $M$ and $\Tr M$ are semi-Gorenstein-projective,
if and only if $M$ is reflexive and both $M$ and $M^*$ are semi-Gorenstein projective:
See (1$'$), (2$'$) and (3).
	\medskip
{\bf Remark 2. Symmetry.} The $\mho$-quiver shows
a symmetry between the semi-Gorenstein-projective
modules and the $\infty$-torsionfree modules:  An indecomposable non-projective module $M$ is
semi-Gorenstein-projective provided there is an infinite $\mho$-path
starting in $M$; and $M$ is $\infty$-torsionfree, provided
there is an infinite $\mho$-path ending in $M$.
	\medskip
{\bf Remark 3. Weakly Gorenstein algebras.} {\it An artin algebra $A$ is left weakly Gorenstein
if and only if there are no modules of $\mho$-type $-\Bbb N,$} see (1$''$).
Similarly, {\it $A$ is right weakly Gorenstein
if and only if there are no modules of $\mho$-type $\Bbb N,$} see (2$''$) and (5).
	\bigskip

{\bf 1.6.} The first example of a semi-Gorenstein-projective module
which is not Gorenstein-projective was constructed by Jorgensen and \c{S}ega [JS] in 2006,
for a commutative algebra of dimension 8. Recently,
Marczinzik [M2] presented some non-commutative algebras
with semi-Gorenstein-projective modules
which are not Gorenstein-projective. In 6.1, we
exhibit a class of 6-dimensional $k$-algebras $\Lambda(q)$ with
parameter $q\in k\setminus\{0\}$ and a family $M(\alpha)$ of 3-dimensional indecomposable $\Lambda(q)$-modules (with $\alpha\in k)$ in order to find new examples:
	\smallskip
{\bf Theorem.} {\it Let $\Lambda(q)$ be the $6$-dimensional algebra defined in $6.1.$
If the multiplicative order of $q$ is infinite,
then the $\Lambda$-modules $M(q)$ and $M(q)^*$ both are
semi-Gorenstein-projective, but $M(q)$ is not torsionless,
thus not Gorenstein-projective;
all the syzygy modules $\Omega^tM(q)$ and $\Omega^t (M(q)^*)$ with $t \ge 0$ are
$3$-dimensional and indecomposable; the module $M(q)^{**}\simeq \Omega M(1)$ is also $3$-dimensional,
but decomposable.}
	\smallskip
{\bf Addendum.} {\it For any $q$, the $\Lambda(q)$-modules
$M(\alpha)$ with $\alpha\in k\setminus q^{\Bbb Z}$
are Gorenstein-projective.
Thus, if $k$ is infinite, then there are infinitely many isomorphism classes of $3$-dimensional
Gorenstein-projective modules.}
	\bigskip
{\bf 1.7. Independence of the total reflexivity conditions.}
It was asked by Avramov and Martsinkowsky [AM] whether the
following conditions which characterize the Gorenstein-projective modules, are independent.
	\smallskip
\item{(G1)} The $A$-module $M$ is semi-Gorenstein-projective.
\item{(G2)} The $A$-dual $M^* = \Hom(M, \ _AA)$ of $M$
is semi-Gorenstein-projective.
\item{(G3)} The $A$-module $M$ is reflexive.
	\medskip
{\bf Theorem.} {\it For artin algebras,
the conditions $(\G1),\ (\G2)$ and $(\G3)$ are independent.}
	\smallskip
Proof. Theorem 1.6 provides a $\Lambda(q)$-module $M$ (namely  $M = M(q)$)
satisfying the conditions (G1), (G2) and not (G3). It remains to use the following
proposition. $\s$
	\medskip
{\bf Proposition.} {\it If a module $M$ is semi-Gorenstein-projective and not Gorenstein-projective,
then $\Omega^2 M$ satisfies {\rm(G1)}  and {\rm(G3)}, but not {\rm(G2)}.

If a module $M'$ satisfies {\rm(G1)} and {\rm(G3)}, but not {\rm(G2)}, then $N = (M')^*$ is a right module satisfying {\rm(G2)} and {\rm(G3)}, but not {\rm(G1)}.}
	\smallskip
Proof.  Let $M$ be semi-Gorenstein-projective.
Then $\Omega^2 M$ is reflexive and semi-Gorenstein-projective.
By Lemma 2.5, $(\Omega^2 M)^* = \Tr M$. Thus $(\Omega^2 M)^*$ is not
semi-Gorenstein-projective (otherwise, $M$ is Gorenstein-projective).

If $M'$ satisfies {\rm(G1)} and {\rm(G3)}, but not {\rm(G2)}, then $(M')^*$ is
reflexive and $(M')^{**} = M'$  is semi-Gorenstein-projective, i.e., $N= (M')^*$ satisfies {\rm(G2)} and {\rm(G3)}, but not {\rm(G1)}.
$\s$
	\medskip
Actually, for our example $A = \Lambda(q)$,
there is also an $A$-module which satisfies (G2), (G3), but not (G1),
namely the module $M(1)$, see 7.5.
	\smallskip
In [JS], where Jorgensen and \c Sega present the first example of a semi-Gorenstein-projective module which is not Gorenstein-projective, they also exhibited
modules satisfying (G1), (G3), but not (G2), and modules satisfying (G2), (G3), but
not (G1). The algebra $A$ considered in [JS] is commutative.
It is an open problem whether there
exists a {\bf commutative} ring $A$ with a module $M$ satisfying (G1), (G2), but not (G3).
	\medskip\smallskip

{\bf 1.8. Outline of the paper.}
The proofs of theorems 1.2, 1.3 and 1.4 are given in sections 2, 3 and 5, respectively.
We use  what we call (as a shorthand)
{\it approximation sequences}, namely  exact
sequences $0 \to X \to Y \to Z \to 0$ with $Y$ projective and $\Ext^1(Z,A) = 0,$
see section 2. Of special interest are the approximation sequences with both $X$ and $Z$
indecomposable and non-projective; in this case, we have
$X = \Omega Z$ and $Z = \mho X$, and we call them {\it $\mho$-sequences,} see section 3.

Section 4 deals with the $\mho$-quiver of $A$.
An essential ingredient in this setting seems to be Corollary 4.4. The corresponding
Remark 1 in 4.4 asserts that the kernel of the canonical map
$\mho^{t}M\to (\mho^{t}M)^{**}$ is equal to $\Ext^{t+1}(\Tr M,A_A)$, and
its cokernel is equal to $\Ext^{i+2}(\Tr M,A_A)$, for all $t\ge 0.$
		
In sections 6 and 7, we present the
6-dimensional algebra $\Lambda = \Lambda(q)$ depending on a parameter $q\in k\setminus\{0\}$,
which we need for Theorem 1.6. We analyze
some 3-dimensional representations which we
denote by $M(\alpha)$ with $\alpha\in k$. Essential
properties of the modules $M(\alpha)$ can be found in 6.3 to 6.5; they are labeled
by (1) to (9). The properties (1) to (5) in 6.3 are those which are needed in order to
exhibit a module, namely $M(q)$,
which is semi-Gorenstein-projective, but not torsionless,
provided the multiplicative order of $q$ is infinite (see 6.4).
The remaining properties (6) to (9) in 6.5
show, in particular, that in case the multiplicative order of $q$ is infinite,
also
the $\Lambda$-dual $M(q)^*$ of $M(q)$
is semi-Gorenstein-projective. The proof of Theorem 1.6 and its Addendum is given
in 6.7 and 6.8.
In 7.1 and 7.2, some components of the
$\mho$-quivers of the algebras $\Lambda$ and $\Lambda^{\op}$ are described.
	
The final sections 8 and 9 add remarks and mention some open questions.
	\medskip
{\bf 1.9. Terminology.}
We end the introduction with some remarks concerning the terminology and its history.
The usual reference for the introduction of Gorenstein-projective modules
(under the name {\it modules of Gorenstein dimension zero})
is the Memoirs by Auslander and Bridger [AB] in 1969. Actually,  in his thesis [Br], Bridger attributes the concept of the Gorenstein dimension to Auslander: In January 1967,
Auslander gave four lectures at the S\'eminaire Pierre Samuel (see the notes [A1] by Mangeney, Peskine and Szpiro). In these lectures, he discussed the
class of all reflexive modules $M$ such that both $M$ and $M^*$ are
semi-Gorenstein-projective modules and denoted it by $G(A)$
([A1, Definition 3.2.2]). Thus $G(A)$ is the class
of the Gorenstein-projective modules and the conditions (G1), (G2) and (G3) served as the first definition.
In [AB, Proposition 3.8], it is shown that a module
$M$ belongs to $G(A)$ if and only if both $M$ and $\Tr M$ are semi-Gorenstein-projective. Of course, we should stress the following:
Whereas some formulations in [AB] assume that the ring $A$
in question is a commutative noetherian ring,
all the essential considerations in [A1, Br, AB] are shown
in the setting of an abelian category with enough projectives,
and of the category of finitely generated modules over a,
not necessarily commutative, noetherian ring.
Enochs and Jenda [EJ1, EJ2] have
reformulated the definition
of Gorenstein-projective modules in terms of complete projective resolutions, see
also [Chr]. Several other names for the Gorenstein-projective modules are in use, they are also
called ``totally reflexive" modules [AM], and ``maximal Cohen-Macaulay" modules [Buch] and ``Cohen-Macaulay" modules [B2].

We should apologize that we propose a new name
for the modules $M$ with $\Ext^i(M,A) = 0$ for all $i\ge 1$, namely {\it semi-Gorenstein-projective}.
These modules have been called for example
``Cohen-Macaulay modules'' or ``stable modules''. However,
the name ``Cohen-Macaulay module'' is in conflict with its established
use for commutative rings, and, in our opinion, the wording ``stable'' may be too vague as
a proper identifier. We hope that the name {\it semi-Gorenstein-projective}
describes well what is going on: that there is something like a half of a
complete projective resolution (``semi'' means ``half'').
We also propose the name {\it left weakly  Gorenstein} for an algebra $A$ with $\gp(A) = {}^\perp A$
(in contrast to ``nearly Gorenstein'' in [M2]);
of course, a Gorenstein algebra $A$
satisfies $\gp(A) = {}^\perp A$, but the algebras with $\gp(A) = {}^\perp A$ seem
to be quite far away from being Gorenstein. The left weakly Gorenstein algebras have also been
called ``algebras with condition (GC)'' in [CH].
	\medskip
{\bf Acknowledgment}. We thank Alex Martsinkovsky for providing copies of [Br] and [A1].
We are indebted to Lars Christensen, Henrik Holm, Zhaoyong Huang,
Rene Marczinzik, Deja Wu for helpful comments.
We also are grateful to two referees for carefully reading 
the manuscript and making valuable suggestions. 
	\bigskip
{\bf 2. Approximation sequences. Proof of Theorem 1.2.}
	\medskip
{\bf 2.1. Lemma.} {\it Let
$\epsilon\:  0 \to X @>\omega>> Y @>\pi>> Z \to 0$
be an exact sequence with $Y$ projective. Then the following conditions are equivalent:
\item{\rm(i)} $\omega$ is a left $\add(A)$-approximation.
\item{\rm(ii)} $\Ext^1(Z,A) = 0.$
\item{\rm(iii)} The $A$-dual sequence $\epsilon^*$ of $\epsilon$ is exact. \par}
	\smallskip
An exact sequence $0 \to X @>>> Y @>>> Z \to 0$ with $Y$ projective satisfying the equivalent properties
will be called in this paper an {\it approximation sequence} (this is just a shorthand,
since it is too vague to be used in general). One may observe that the conditions (i), (ii) and (iii)
are equivalent for any exact sequence $\epsilon\: X @>\omega>> Y @>>> Z \to 0$ with $Y$ projective,
even if
$\omega$ is not injective, but we are only interested in the short exact sequences.
	\smallskip
Proof of the equivalence of the properties.
Since $Y$ is projective, applying $\Hom(-,A)$ to $\epsilon$ we get the exact sequence
$0 \to Z^* @>\pi^*>> Y^* @>\omega^*>> X^* @>>> \Ext^1(Z,A) \to 0$. Note that $\omega$ is
a left $\add(A)$-approximation if and only if $\omega^*$ is surjective. From this we get the
equivalence of (i) and (ii) and the equivalence of (ii) and (iii).
\hfill$\square$
	\medskip
{\bf 2.2.} Also the following basic lemma is well-known (see, for example [R]).
	\medskip
{\bf Lemma.} {\it Let $P_{-1} @>f>> P_0 @>g>> P_1$ be an exact sequence of
projective modules and let $g = up$ be a factorization with $p:P_0 \to I$ epi and
$u\:I \to P_1$ mono. Then
$P_{-1}^* @<f^*<< P_0^* @<g^*<< P_{1}^*$ is exact if and only if $u$ is a left $\add(A)$-approximation.}
	\medskip
For the convenience of the reader, we insert the proof. Since $f^*g^* = (gf)^* = 0,$
we have $\Im g^* \subseteq \Ker f^*.$ Assume now that $u$ is a left $\add(A)$-approximation
and let $h\in \Ker f^*,$ thus $hf = 0.$ Since $p$ is a cokernel of $f$, there is
$h'$ with $h = h'p.$ Since $u$ is a left $\add(A)$-approximation, there is $h''$ with
$h' = h''u.$ Thus $h = h'p = h''up = h''g = g^*(h'')$ belongs to the image of $g^*,$
there also $\Ker f^* \subseteq \Im g^*$.

Conversely, we assume that $\Im g^* = \Ker f^*$ and let $h\:I \to A$ be a map.
Then $hpf = 0$, so that $f^*(hp) = 0$. Therefore $hp$ belongs to $\Ker f^*$, thus
to $\Im g^*$. There is $h''\in P_1^*$ with $hp = g^*(h'') = h''g = h''up,$
and therefore $h = h''u.$
\hfill$\square$
	\medskip
This Lemma will be used in various settings, see 4.3.
	\medskip
{\bf 2.3.} {\it A semi-Gorenstein-projective and $\Omega$-periodic module is Gorenstein-projective.}
	\medskip
Proof. Let $M$ be semi-Gorenstein-projective and assume that $\Omega^tM = M$ for some
$t\ge 1$. Let $\cdots \to P_{i} \to \cdots \to P_0 \to M \to 0$  be a minimal projective
resolution of $M$. Then
$$
 0 \to \Omega^tM \to P_{t-1} \to \cdots \to P_0 \to M \to 0 \tag{+}
$$
is the concatenation of approximation sequences. Since $\Omega^tM = M$, we can concatenate
countably many copies of $(+)$ in order to obtain a double infinite acyclic chain complex of projective
modules. As a concatenation of approximation sequences, it is a complete projective
resolution. Therefore, $M$ is Gorenstein projective.
\hfill$\square$
	\medskip
{\bf 2.4.} Here are two essential observations.
	\medskip
{\bf (a)} {\it Let $0 \to X \to Y \to Z \to 0$ be an approximation sequence. Then
$\phi_X$ is surjective if and only if $Z$ is torsionless.} We can also say: {\it $X$ is reflexive if
and only if $Z$ is torsionless.}
	\smallskip
{\bf (b)} {\it Let $0 \to X \to Y \to Z \to 0$ be an approximation sequence. Then
$\Ext^1(X^*,A_A) = 0$ if and only if $\phi_Z$ is surjective.}
	\medskip
Proof of (a) and (b). Since $0 \to X @>\omega>> Y @>\pi>>  Z \to 0$ is an approximation sequence,
it follows that
$$
 0\longrightarrow Z^* \overset {\pi^*}\to \longrightarrow Y^* \longrightarrow X^* \longrightarrow 0
$$
is an exact sequence of right $A$-modules.
This induces an exact sequence
$$
 0\longrightarrow X^{**} \longrightarrow Y^{**}
\overset {\pi^{**}} \to \longrightarrow Z^{**} \longrightarrow
\Ext_A^1(X^*, A_A) \longrightarrow 0
$$
of left $A$-modules,
and we obtain the commutative diagram
$$\CD
  0 @>>> X @>>> Y @>\pi>> Z @>>> 0 \cr
  @.    @VV\phi_XV   @|        @VV\phi_ZV  \cr
  0 @>>> X^{**} @>>> Y^{**} @>\pi^{**}>> Z^{**} @>>> \Ext^1(X^*,A_A) @>>> 0\ .
\endCD
$$
By the Snake Lemma, the kernel of $\phi_Z$ is isomorphic to the cokernel of $\phi_X$,
Thus $\phi_Z$ is a monomorphism if and only if $\phi_X$ is an epimorphism.  Since $X$ is torsionless, $X$ is reflexive if and only if $\phi_X$ is surjective. This is (a).

By the commutative diagram above, we see that $\phi_Z$ is epic if and only if so is $\pi^{**}$,
and if and only if $\Ext^1_A(X^*,A_A) = 0$. This is (b).
\hfill$\square$
	\medskip
{\bf Corollary.} {\it A module $X$ is reflexive if and only if both $X$ and $\mho X$ are torsionless.}
	\medskip
Proof. If $X$ is reflexive, then it is torsionless. Thus we may assume from the
beginning that $X$ is torsionless. Any minimal left
$\add(A)$-approximation $X \to Y$ is injective and its cokernel is $\mho X$.
The exact
sequence $0\to X \to Y \to \mho X \to 0$ is an approximation sequence, and
2.4 (a) asserts that $X$ is reflexive iff $\mho X$ is torsionless.
$\s$
	\medskip
{\bf Remark.} The assertion of the corollary can be strengthened as follows.
For any module $X$, let us denote by $\K X$ the kernel of the map
$\phi_X\:X \to X^{**}$. Of course,  $\K X$ is the kernel of any left
$\add(A)$-approximation of $X$. Therefore
$X$ is torsionless if and only if $\K X = 0.$ Claim:
{\it The cokernel of the map $\phi_X\:X \to X^{**}$ is isomorphic to $\K\mho X.$}
	\medskip
Proof. Let $u\: X \to Y$ be a minimal $\add(A)$ approximation, say with cokernel
$p\:Y \to \mho X$. The $A$-dual of the exact sequence $X @>u>> Y @>p>> \mho X @>>> 0$
is $0 @<<< X^* @<u^*<< Y^* @<p^*<< (\mho X)^* @<<< 0,$ since $u$ is an $\add(A)$-approximation.
Using again $A$-duality, we obtain the exact sequence
$0 @>>>X^{**} @>u^{**}>> Y^{**} @>p^{**}>> (\mho X)^{**}$. Thus there is the following
commutative diagram with exact rows:
$$
\CD
    @.        X @>u>> Y @>\pi>> \mho X @>>> 0 \cr
  @.    @VV\phi_X V   @VV\phi_Y V   @VV\phi_{\mho X}V  \cr
  0 @>>> X^{**} @>u^{**}>> Y^{**} @>\pi^{**}>> {\mho X}^{**}.
\endCD
$$
Since $\phi_Y$ is an isomorphism, the snake lemma yields $\Cok\phi_X \simeq
\Ker(\phi_{\mho X}) = \K \mho X$. $\s$
	\medskip
In 4.4, we will rewrite both $\K X$ and $\K\mho X$ in order to obtain the
classical Auslander-Bridger sequence (see Corollary and Remark 1 in 4.4).
	\medskip
{\bf 2.5. Lemma.} {\it Let $M$ be a module with $\Ext^i(M,A) = 0$ for $i=1,2$.
Then $\Tr M \simeq (\Omega^2 M)^*$ and there is a projective module $Y$ such that
$M^* \simeq \ \Omega^2\Tr M \oplus Y.$}
	\medskip
Proof: Let $\pi\:PM \to M$ and $\pi'\:P\Omega M \to \Omega M$ be projective covers with inclusion maps
$\omega\:\Omega M \to PM$ and
$\omega'\:\Omega^2 M \to P\Omega M$.
Then $\omega\pi'$ is a minimal projective
presentation of $M$. By definition, $\Tr M$ is the cokernel of $(\omega\pi')^*$.
Since $\Ext^i(M,A) = 0$ for $i=1,2$, the exact sequences
$$
 0 @>>> \Omega^2M @>\omega'>> P\Omega M @>\pi'>> \Omega M @>>> 0,\quad
 0 @>>> \Omega M @>\omega>> P M @>\pi>>  M @>>> 0
$$
are approximation sequences. As a consequence, the corresponding $A$-dual sequences
$$
 0 @<<< (\Omega^2M)^* @<(\omega')^*<< (P\Omega M)^* @<(\pi')^*<< (\Omega M)^* @<<< 0,\quad
 0 @<<< (\Omega M)^* @<\omega^*<< (P M)^* @<\pi^*<<  M^* @<<< 0
$$
are exact. The concatenation
$$ \epsilon\:\qquad
 0 @<<< (\Omega^2M)^* @<(\omega')^*<< (P\Omega M)^* @<(\omega\pi')^*<< (P M)^* @<\pi^*<<  M^* @<<< 0
$$
shows that $(\Omega^2M)^*$ is a cokernel of $(\omega\pi')^*,$ thus
$\Tr M \simeq (\Omega^2M)^*$. In addition, $\epsilon$ shows that $\Omega^2 \Tr M
= \Omega^2(\Omega^2M)^*$
is the direct sum of $M^*$ and a projective module $Y$.
\hfill$\square$

	\medskip
{\bf 2.6.} Proof of Theorem 1.2.

(1) implies (2) to (7): This follows directly from well-known properties of Gorenstein-projective
modules. Namely, assume (1) and
let $M$ be Gorenstein-projective. Then $M$ is reflexive, this yields (3), but, of course,
also (2) and (4). Second,
$M^*$ is Gorenstein-projective, thus semi-Gorenstein-projective, therefore we get (5) and (6).
Finally, $\Tr M$ is Gorenstein-projective, thus semi-Gorenstein-projective, therefore we get (7).
	
Both (3) and (4) imply (2): Let $M$ be semi-Gorenstein-projective. Consider the approximation sequence
$0 \to \Omega M \to PM \to M \to 0$ and note that $\Omega M$ is again semi-Gorenstein-projective.
If (3) or just (4) holds, we know that $\phi_{\Omega M}$ is surjective, thus by 2.4 (a), $M$ is
torsionless.
	
Both (6) and (7) imply (2):
Let $M$ be semi-Gorenstein-projective. Consider the approximation sequences
$0 \to \Omega M \to PM \to M \to 0$ and
$0 \to \Omega^2 M \to PM \to \Omega M \to 0$. Since $M$ is semi-Gorenstein-projective,
also $\Omega^2 M$ is semi-Gorenstein-projective.
If (6) holds, we use (6) for $\Omega^2 M$ in order to see that $\Ext^1((\Omega^2M)^*,A_A) = 0$.
If (7) holds, we use (7) for $M$ in order to see that $\Ext^1(\Tr M,A_A) = 0.$
According to 2.5, we see that $\Tr M = (\Omega^2 M)^*.$
Thus in both cases (6) and (7), we have $\Ext^1((\Omega^2M)^*,A_A) = 0$. According to 2.4 (b), it follows from
$\Ext^1((\Omega^2M)^*,A_A) = 0$ that $\phi_{\Omega M}$ is surjective. By 2.4 (a), $M$ is torsionless.
	
Trivially, (5) implies (6). Altogether we have shown that any one of the
assertions (3) to (7) implies (2).
	
It remains to show that (2) implies (1). Let $M$ be semi-Gorenstein-projective and torsionless.
We want to show that $M$ is Gorenstein-projective.
Let $M_i =  \mho^iM$ for all $i\ge 0$ (with $M_0 = M$). Since $M_0$ is torsionless,
there is an  approximation sequence $0 \to M_0 \to P_1 \to M_1 \to 0$, and $M_1$ is again
semi-Gorenstein-projective. By assumption, $M_1$ is again torsionless.
Inductively, starting with a torsionless module $M_i$,
we obtain an approximation sequence $\epsilon_i\: 0 \to M_i \to P_{i+1} \to M_{i+1} \to 0$,
we conclude that with $M_i$ also $M_{i+1}$ is semi-Gorenstein-projective. By
(2) we see that $M_{i+1}$ is torsionless, again.
Concatenating a minimal projective resolution of $M$ with these
approximation sequences $\epsilon_i$, for $0\le i$,
we obtain a complete projective resolution of $M$.
\hfill$\square$
	\bigskip
{\bf 3. $\mho$-sequences. Proof of theorem 1.3.}
	\medskip
{\bf 3.1.} An approximation sequence $0 \to X \to Y \to Z \to 0$ will be called an
{\it $\mho$-sequence} provided
that both $X$ and $Z$ are indecomposable and not projective (the relevance of such
sequences was stressed already in [RX]).
	\medskip
{\bf Lemma.} {\it An approximation sequence is the direct sum of $\mho$-sequences and
an exact sequence $0\to X'\to Y'\to Z'\to 0$ with $X',\ Z'$ (thus also $Y'$) being projective.}
	\medskip
Proof: Let $0 \to X @>\omega>> Y @>\pi>> Z \to 0$ be an approximation sequence.
Since $Y$ is projective and $\pi$ is surjective, a direct decomposition
$Z = Z_1\oplus Z_2$ yields a direct sum decomposition of the sequence.
Since $\omega$ is a left $\add(A)$-approximation, there is also the corresponding assertion:
If $X = X_1\oplus X_2$, then $X @>\omega>> Y$ is the direct sum of two maps $X_1\to Y_1$ and $X_2 \to Y_2$,
thus again we obtain a direct sum decomposition of the sequence.
This shows that for an indecomposable approximation sequence $0 \to X @>\omega>> Y @>\pi>> Z \to 0$,
the modules $X$ and $Z$ are indecomposable or zero (and, of course, not both can be zero).

If $Z$ is indecomposable and projective, then the sequence $0 \to X \to Y \to Z \to 0$ splits off
$0 \to 0 \to Z @>1>> Z \to 0$, thus $X = 0.$ Similarly, if $X$ is indecomposable and projective, then
the sequence $0 \to X \to Y \to Z \to 0$ splits off
$0 \to X @>1>> X \to 0\to 0$, thus $Z = 0.$

It remains that $0 \to X \to Y \to Z \to 0$ is an approximation sequence with both $X$ and $Z$
being indecomposable and non-projective.
\hfill$\square$
	\medskip
{\bf 3.2. Lemma.} {\it Let $\epsilon\:0 \to X @>\omega>> Y @>\pi>> Z \to 0$ be an exact sequence. The following conditions are
equivalent:
\item{\rm (i)}  $\epsilon$ is an $\mho$-sequence.
\item{\rm (ii)}
    $X$ and $Z$ are indecomposable and not projective, $\omega$ is a minimal
    left $\add(A)$-approximation, $\pi$ is a projective cover, $X = \Omega Z$, $Z = \mho X$.
\item{\rm (iii)} $X$ is indecomposable and not projective, $\omega$ is a minimal left $\add(A)$-approximation.
\item{\rm (iv)} $Z$ is indecomposable and not projective, $\pi$ is a projective cover, and $\Ext^1(Z,A) = 0.$
\item{\rm (v)} $X = \Omega Z$, $Y$ is projective, $Z = \mho X$, and $X$ is indecomposable.
\item{\rm (vi)} $X = \Omega Z$, $Y$ is projective, $Z = \mho X$, and $Z$ is indecomposable.\par}
	\medskip
Proof: (i) implies (ii):
Let $\epsilon$ be an $\mho$-sequence. Then $\omega$
has to be minimal, since otherwise $\epsilon$ would split off a non-zero sequence
of the form $0\to 0 \to P @>1>> P \to 0$ with $P$ projective.
Similarly, $\pi$ has to be a projective cover, since
otherwise $\epsilon$ would split off a non-zero sequence
of the form $0 \to P @>1>> P \to 0 \to 0$.
Since $\omega$ is a minimal left $\add(A)$-approximation and $Z$ is the cokernel of
$\omega$, we see that $Z = \mho X$. Since $\pi$ is a projective cover of $Z$ and $X$ is its kernel, $X = \Omega Z$.

(ii) collects all the relevant properties of an $\mho$-sequence.
The condition (iii), (iv), (v) and (vi) single out some of these properties,
thus (ii) implies these conditions.

(iii) implies (i): Since $X$ is indecomposable and not projective, $\epsilon$ has no direct summand
$0\to P @>1>> P \to 0\to 0 $.
Since $\omega$ is left minimal, $\epsilon$ has no direct summand $0\to 0 \to P @>1>> P \to 0$.
Similarly, (iv) implies (i).

Both (v) and (vi) imply (i): Since $Z = \mho X$, we have $\Ext^1(Z,A) = 0.$ This shows that the sequence
is an approximation sequence. Since $X = \Omega Z$, the sequence  $\epsilon$ has no direct summand
of the form $0\to P @>1>> P \to 0\to 0 $.
Since $Z = \mho X$, the sequence  $\epsilon$ has no direct summand of the form
$0\to 0 \to P @>1>> P \to 0$. Thus, $\epsilon$ is a direct sum of $\mho$-sequences.
Finally, since $X$ or $Z$ is indecomposable, $\epsilon$ is an $\mho$-sequence.
\hfill$\square$

	\medskip
{\bf 3.3. Corollary.}
{\it If $M$ is indecomposable, non-projective, semi-Gorenstein-projective, then
$\Omega M$ is
indecomposable, non-projective, semi-Gorenstein-projective and $M = \mho\Omega M.$}
	\medskip
Proof. Since $M$ is semi-Gorenstein-projective module, the canonical sequence
$\epsilon\:0 \to \Omega M \to P M \to M \to 0$
is an approximation sequence. Since $M$ is indecomposable and not projective, and $PM \to M$ is
a projective cover, $\epsilon$ is an $\mho$-sequence, thus $\Omega M$ is indecomposable and non-projective,
and $M = \mho\Omega M$, by 3.2. Of course, with $M$ also $\Omega M$ is semi-Gorenstein-projective.
\hfill$\square$
	\medskip
{\bf 3.4. Lemma.}
{\it If the number of isomorphism classes of indecomposable modules which are both
semi-Gorenstein-projective and torsionless is finite, then
any indecomposable non-projective semi-Gorenstein-projective module $M$ is $\Omega$-periodic.}
	\medskip
Proof. According to 3.3,
the modules $\Omega^t M$ with $t\ge 1$ are indecomposable modules which are torsionless
and semi-Gorenstein-projective.
Since there are only finitely many isomorphism classes of indecomposable torsionless
semi-Gorenstein-projective modules, there are natural numbers $1 \le s < t$ with
$\Omega^s M = \Omega^t M$. Then $M = \mho^s\Omega^s M = \mho^s \Omega^t M = \Omega^{t-s}M$
and $t-s \ge 1$, thus $M$ is $\Omega$-periodic.
\hfill$\square$
	\medskip
{\bf 3.5. Proof of Theorem 1.3.} We assume that the number of isomorphism classes of indecomposable
torsionless semi-Gorenstein-projective modules is finite. According to 3.4,
any indecomposable non-projective semi-Gorenstein-projective module is $\Omega$-periodic.
2.3 shows that any semi-Gorenstein-projective $\Omega$-periodic
module is Gorenstein-projective.
\hfill$\square$
	\medskip
{\bf Remark.} One of the referees has pointed out that
Theorem 1.3 can be improved by replacing the class
of all torsionless modules by an arbitrary full subcategory $\Cal X$ which is closed under direct
summands, contains $\add(A)$, and contains for any indecomposable module $M$ at least one
syzygy module $\Omega^n M$. {\it If ${}^\perp A \cap \Cal X$ contains only finitely
many isomorphism classes of indecomposable modules,
then $\Lambda$ is left weakly Gorenstein and any Gorenstein-projective module
is $\Omega$-periodic.}

	\medskip
{\bf 3.6. Torsionless-finite algebras.}
An artin algebra $A$ is said to be {\it torsionless-finite} if there are only finitely many
isomorphism classes of indecomposable torsionless modules. Theorem 1.3 implies
that
{\it any torsionless-finite artin algebra is left weakly  Gorenstein,}  as Marczinzik [M1] has
shown.
If $\Lambda$ is torsionless-finite, also $\Lambda^{\op}$ is torsionless-finite [R],
thus a torsionless-finite artin algebras is also right weakly Gorenstein.
Note that many interesting classes of algebras are known to be torsionless-finite.
In particular, we have
	\smallskip
{\it The following algebras are torsionless-finite, hence left  and right weakly  Gorenstein.
	
\item{\rm(1)} Algebras $A$ such that $A/\soc (_AA)$ is representation-finite.

\item{\rm(2)} Algebras stably equivalent to hereditary algebras, in particular all
 algebras with radical square zero.

\item{\rm(3)} Minimal representation-infinite algebras.

\item{\rm(4)} Special biserial algebras without indecomposable projective-injective modules.\par}
	\medskip
See for example [R], where also other classes of torsionless-finite algebras are listed.
	\medskip
Chen [Che] has shown that a connected algebra $A$ with radical square zero either is
self-injective, or else all the Gorenstein-projective modules are projective.
The assertion
that algebras with radical square zero are weakly  Gorenstein complements this result.
	
\vfill\eject
{\bf 4. The $\mho$-quiver.}
	\medskip
{\bf 4.1.} We recall that the $\mho$-quiver of $A$ has as vertices the isomorphism classes $[X]$ of the
indecomposable non-projective modules $X$ and {\it there is an  arrow
$$
{\beginpicture
\setcoordinatesystem units <1cm,1cm>
\put{$[X]$} at 0 0
\setdashes <1mm>
\arr{1.5 0}{0.5 0}
\put{$[Z]$} at 2 0
\endpicture}
$$
provided that $X$ is torsionless and $Z = \mho X,$ thus provided that
there exists an $\mho$-sequence
$0 \to X \to Y \to Z \to 0.$} We will also write the vertex $[X]$ simply as $X$.
	\medskip
{\bf 4.2. The $A$-dual of an $\mho$-sequence.}
	\medskip
{\bf Lemma.} (a) {\it Let $\epsilon\:0 \to X \to Y \to Z \to 0$
be an approximation sequence  and assume that
$X$ is reflexive. Then $\Ext^1(X^*,A_A) = 0$ if and only if $Z$ is reflexive,
if and only if
the $A$-dual $\epsilon^*$ of $\epsilon$ is again an approximation sequence.}
	\smallskip
(b) {\it Let $\epsilon\:0 \to X \to Y \to Z \to 0$
be an $\mho$-sequence with $X$ reflexive. Then $Z$ is reflexive,
if and only if the $A$-dual $\epsilon^*$ of $\epsilon$ is again an $\mho$-sequence.}
	\medskip
Proof. (a)
By 2.4 (a), we see that $Z$ is always torsionless. Thus 2.4 (b) shows that
$\Ext^1(X^*,A_A) = 0$ if and only if $Z$ is reflexive. First,
assume that $Z$ is reflexive. Then $\Ext^1(X^*,A_A) = 0$, and therefore
we see that the $A$-dual sequence $\epsilon^*$ is exact.
We dualize a second time: the sequence $\epsilon^{**}$ is isomorphic to the sequence $\epsilon,$
since the three modules $X,Y,Z$ are reflexive. This means that $\epsilon^{**}$ is exact, and therefore $\epsilon^*$ is an approximation sequence.
Second, conversely, if $\epsilon^*$ is an approximation sequence, then it is
exact, and therefore $\Ext^1(X^*,A_A) = 0$, thus $Z$ is reflexive.

(b) Assume now that $\epsilon$ is an $\mho$-sequence. First, assume
that $Z$ is reflexive. Since $X,Z$ both are reflexive,
indecomposable and non-projective, also $X^*$ and $Z^*$ are indecomposable and non-projective, as we will show below. Thus $\epsilon^*$ is an $\mho$-sequence.
Conversely, if $\epsilon$ is an $\mho$-sequence, then it
is an approximation sequence and thus $Z$ is reflexive by (a).
\hfill$\s$
	\medskip
We have used some basic facts about the $A$-dual $M^*$ of a module $M$.
	\smallskip
\item{(1)} {\it $M^*$ is always torsionless.}

\item{(2)} {\it If $M$ is non-zero and torsionless, then $M^*$ is non-zero.}

\item{(3)} {\it If $M$ is reflexive, indecomposable and non-projective, then $M^*$ is reflexive,
indecomposable and non-projective.\par}
	\smallskip

Here are the proofs (or see for example [L,  p.144]).
(1) There is a surjective map $p\:P \to M$ with $P$ projective. Then $p^*\:M^* \to P^*$
is an embedding of $M^*$ into the projective module $P^*$. The assertion (2) is obvious.

(3) Let $M$ be reflexive, indecomposable and non-projective. Consider a direct decomposition
$M^* = N_1\oplus N_2$ with $N_1\neq 0$ and $N_2 \neq 0$. Since $M^*$ is torsionless by (1),
both modules $N_1$ and $N_2$ are torsionless, therefore $N_1^* \neq 0, N_2^* \neq 0,$
thus there is a proper direct decomposition
$M^{**} = N_1^* \oplus N_2^*$. Since $M$ is reflexive and indecomposable, this is impossible.
Thus $M^*$ has to be
indecomposable. If $M^*$ is projective, then also $M^{**}$ is projective. Again, since $M$ is
reflexive, this is impossible.

It remains to show that $M^*$ is reflexive. Since $M^{**}$ is isomorphic to $M$, we see
that $M^{***}$ is isomorphic to $M^*$, thus the canonical map $M^* \to M^{***}$
has to be an isomorphism (since it is a monomorphism of modules of equal length).
\hfill$\square$
\medskip

{\bf 4.3.}
Lemma 2.2 outlines the importance of left $\add(A)$-approximations when dealing
with exact sequences of projective modules. Let us  give a unified treatment
of the relevance of approximation sequences and of $\mho$-sequences.
	\medskip
(a) {\it An exact sequence
$\cdots \longrightarrow P^{-1}\longrightarrow P^{0}
  \overset {d^0} \to {\longrightarrow} P^{1}\longrightarrow \cdots$
is a complete projective resolution if and only
if it is the concatenation of approximation sequences.}
	\smallskip
(a$'$) {\it An indecomposable non-projective module $M$ is Gorenstein-projective
if and only if $[M]$ is the start of an infinite $\mho$-path and the end of an infinite
$\mho$-path.}
 	\smallskip
(b) {\it A module $M$ is semi-Gorenstein-projective if and only if a projective
resolution (or, equivalently, any projective resolution) is the concatenation
of approximation sequences.}
	\smallskip
(b$'$) {\it An indecomposable non-projective module $M$ is semi-Gorenstein-projective
if and only if $[M]$ is the start of an infinite $\mho$-path.}
	\smallskip
(c) {\it A module $M$ is
reflexive and $M^*$ is semi-Gorenstein-projective if and only if
there is an exact sequence $0 \to M \to P^1 \to P^2 \to \cdots$
which is the concatenation of approximation sequences.}
	\smallskip
(c$'$) {\it An indecomposable non-projective module $M$ is
reflexive and $M^*$ is semi-Gorenstein-projective
if and only if $[M]$ is the end of an infinite $\mho$-path.}
	\medskip
Proof: We use that the
$A$-dual of an approximation sequence is exact, thus
the $A$-dual of the concatenation of approximation sequences is exact.
	\medskip
(a) Let $P^\bullet$ be a double infinite exact sequence of projective modules
with maps $d^i\:P^i \to P^{i+1}$
Write $d^i = \omega^{i}\pi^{i}$ with $\pi^{i}$ epi
and $\omega^{i}$ mono. If $P^\bullet$ is a complete projective resolution,
then the exactness of $(P^\bullet)^*$ at $(P^i)^*$
implies that $\omega^{i}$ is a left $\add(A)$-approximation, see 2.2.
Thus $P^{\bullet}$ is the concatenation of approximation sequences.
	\smallskip
(b) Let $\cdots \to P_i \to \cdots \to P_1 \to P_0 \to M \to 0$ be a projective
resolution of $M$.
write the map $P_{i+1} \to P_i$ as $\omega_i\pi_i$ with
$\pi_{i}$ epi and $\omega_{i}$ mono. If the $A$-dual of the sequence
$\cdots \to P_i \to \cdots \to P_0$ is exact, then all the maps $\omega_i$ with $i\ge 1$ have to be left $\add(A)$-approximations.
This shows that the projective resolution
is the concatenation of approximation sequences.
	\smallskip
(b$'$) Let $M$ be indecomposable, non-projective and semi-Gorenstein-projective.
Since $\Ext^1(M,A) = 0$, the sequence $0 \to \Omega M \to PM \to M \to 0$
is an $\mho$-sequence and $\Omega M$ is again indecomposable and non-projective.
Also, $\Omega M$ is semi-Gorenstein-projective. Thus, we can iterate the procedure
and obtain the infinite path
$$
{\beginpicture
\setcoordinatesystem units <1cm,1cm>
\put{$[M]$} at 0 0
\setdashes <1mm>
\arr{-.5 0}{-1.3 0}
\put{$[\Omega M]$} at -2 0
\arr{-2.7 0}{-3.3 0}
\put{$[\Omega^2 M]$} at -4 0
\arr{-4.7 0}{-5.3 0}
\put{$\cdots$} at -6 0
\put{$(*)$} at -8 0
\endpicture}
$$
Conversely, assume that there is an infinite path starting with $[M]$, then it is
of the form $(*)$.
Thus, for all $i\ge 1$, we have  $\Ext^i(M,A) \simeq \Ext^1(\Omega^{i-1}M,A) = 0.$
	\medskip
Proof of (c) and (c$'$).  Assume that there are given approximation sequences
$\epsilon_i\: 0 \to M^i \to P^{i+1} \to M^{i+1} \to 0$
for all $i\ge 0$, with $M^0 = M.$
Then all the modules $M^i$ are torsionless, thus reflexive by 2.4(a). In particular,
$M$ itself is reflexive.
The $A$-dual of $\epsilon_i$ is the sequence
$$
 \epsilon_i^*\: 0 @<<< (M^i)^* @<<< (P^{i+1})^* @<<< (M^{i+1})^* @<<< 0,
$$
which again is an approximation sequence by 4.2(a).
The concatenation of the sequences $\epsilon_i^*$ is a projective resolution of
$M^* = (M^0)^*$. According to (b), $M^*$ is
semi-Gorenstein-projective, since all the sequences $\epsilon_i^*$ are
approximation sequences.

\smallskip

Conversely, assume that $M$ is reflexive and $M^*$ is
semi-Gorenstein-projective. We want to construct a sequence
$0 \to M \to P^1 \to P^2 \to \cdots$ which is the
concatenation of approximation sequences. It is sufficient to consider
the case where $M$ is indecomposable (in general, take the direct
sum of the sequences).
If $M$ is projective, then $0 \to M \to M \to 0 \to \cdots$ is the
concatenation of approximation sequences.

Thus, it remains to consider the case where $M$ is indecomposable and
not projective.
Since $M$ is torsionless, there is an $\mho$-sequence
$\epsilon_0\:0 \to M \to P^1 \to M^1 \to 0$ (with $M^1 = \mho M$).
Note that $M^1$ is indecomposable,
not projective, and that the $A$-dual
$\epsilon_0^*\:0 @<<< M^* @<<< (P^1)^* @<<< (M^1)^* @<<< 0$ is exact.
Since $M$ is reflexive, $M^1$ is torsionless by 2.4 (a).
Since $M^*$ is semi-Gorenstein-projective, $\Ext^1(M^*,A_A) = 0$, therefore
$\phi_{M^1}$ is surjective and $\epsilon_0^*$ is an $\mho$-sequence, by 4.2.
Altogether we know now that $M^1$ is reflexive, but also that $(M^1)^* =
\Omega(M^*).$
With $M^*$ also $\Omega(M^*)$ is semi-Gorenstein-projective.

Thus $M^1$ satisfies again the assumptions of being
indecomposable, not projective, reflexive and that its $A$-dual $(M^1)^*$
is semi-Gorenstein-projective.
Thus we can iterate the procedure for getting the next $\mho$-sequence
$\epsilon_1\:0 \to M^1 \to P^2 \to M^2 \to 0$, with $M^2 = \mho^2 M$,
and so on. Altogether, we obtain the infinite path:
$$
{\beginpicture
\setcoordinatesystem units <1cm,1cm>
\put{$[M]$} at 0 0
\setdashes <1mm>
\arr{1.3 0}{.5 0}
\put{$[\mho M]$} at 2 0
\arr{3.3 0}{2.7 0}
\put{$[\mho^2 M]$} at 4 0
\arr{5.3 0}{4.7 0}
\put{$\cdots$} at 6 0
\endpicture}
$$
This completes the proof of (c$'$) and thus also of (c).
	\smallskip
(a$'$) This follows immediately from (b$'$) and (c$'$).
\hfill$\square$

\medskip

{\bf 4.4.} For any module $M$, we have denoted by $\K M$ the kernel of $\phi_M\:M \to M^{**}$.
We are going to identify $\K M$ with $\Ext^1(\Tr M, A_A)$. Compare [A2, Proposition 6.3]. As a consequence, we see that  $\mho M = \Tr\Omega\Tr M$.

	\medskip
{\bf Lemma.} {\it Let $M$ be a module. Then $\Ext^1(\Tr M,A_A) \simeq \K M$ and there is
a right module $Q$ such that $\Omega \Tr M \simeq \Tr \mho M\oplus  Q$.
As a consequence, $\mho M \simeq \Tr\Omega\Tr M,$ thus  
$\mho^t(M)\cong \Tr \Omega^t \Tr (M)$ for $t\ge 1.$}
	\smallskip
Proof.
Let $P^{0} @>f>> P^1 @>p>> M @>>> 0$ be a minimal projective presentation of $M$. Thus $\Tr M$ is the cokernel of $f^*.$ Let $g'\:M \to P^2$ be a minimal left
$\add(A)$-approximation. Then $\K M$ is the kernel
of $g'$, thus $g' = uq$, where $q\:M \to M/\K M$ is the canonical projection and $u$ is
injective. Let $g = g'p = uqp.$ The composition
$$
 P^0 @>f>> P^1 @>g>> P^2
$$
is zero and the homology $H(P^0 @>f>> P^1 @>g>> P^2)$ is just $\K M,$ since
$\Ker(g)/\Im(f) \simeq \Ker(qp)/\Ker (p) \simeq \K M.$

We claim that the $A$-dual sequence
$$
  (P^0)^* @<f^*<< (P^1)^* @<g^*<< (P^2)^* \tag{$*$}
$$
is exact. Since $gf = 0$, we have $f^*g^* = 0.$ Conversely,
let $h\:P^1 \to A$ be in the kernel of $f^*,$
thus $hf = 0.$ Therefore $h$ factors through $p = \Cok f$,
say $h = h'p$ with $h'\: M \to A$. Since $uq$ is a left $\add(A)$-approximation,
we obtain $h''\:P^{2} \to A$ with $h' = h''uq$. Thus
$h = h'p = h''uqp = h''g = g^*(h'')$ is in the image of $g^*.$

Since the cokernel of $f^*$ is $\Tr M$, it follows that
$(*)$ is the begin of a projective resolution of $\Tr M$  and hence
$\Ext^1(\Tr M,A_A)$ is obtained by applying $\Hom(-,A)$ to $(*)$
and taking the homology at the position 1.
Applying $\Hom(-,A)$ to $(*)$ we retrieve the sequence $P^0 @>f>> P^1 @>g>> P^2$,
thus $\Ext^1(\Tr M,A_A)$ is equal to $H(P^0 @>f>> P^1 @>g>> P^2) \simeq \K M$.
This is the first assertion.
	\medskip
By definition, the cokernel of $g'$ is $\mho M.$ Thus
the cokernel of $g$ is $\mho M$, and therefore $\Cok g^*\simeq \Tr \mho M \oplus Q'$
for some projective right module $Q'$. Now $\Cok g^* = \Im f^*$, since
$(*)$ is exact. Since $\Cok f^* = \Tr M$, we have $\Omega \Tr M \simeq \Im f^*\oplus Q''$
for some projective right module $Q''$. This shows
that $\Omega \Tr M \simeq \Im f^*\oplus Q'' = \Cok g^*\oplus Q'' \simeq
\Tr \mho M\oplus Q'\oplus Q'' = \Tr\mho M \oplus Q$ with $Q = Q'\oplus Q''$.
This is the second assertion.
	\medskip
Applying $\Tr$ to the isomorphism $\Omega\Tr M \simeq \Tr\mho M \oplus Q$, one obtains
$\Tr\Omega\Tr M \simeq \Tr(\Tr \mho M\oplus Q) = \Tr\Tr \mho M$. Since $\mho M$ has no
non-zero projective direct summand, one gets $\Tr\Tr\mho M \simeq \mho M$.
Thus $\mho M \simeq \Tr\Tr \mho M \simeq \Tr\mho\Tr M.$
$\s$

	\medskip
{\bf Corollary.} {\it Let $M$ be a module. Then for all $t\ge 0$ one has
$$
  \Ext^{t+1}(\Tr M,A_A) \simeq \K(\mho^{t}M).
$$
In particular, $\mho^{t}M$ is torsionless if and only if $\Ext^{t+1}(\Tr M,A_A) = 0.$
Also, $\Omega^t \Tr M \simeq \Tr \mho^t M\oplus Q_t$ for some projective right module $Q_t$.}
	\smallskip
Proof. By induction on $t$, one has $\Omega^t \Tr M \simeq \Tr \mho^t M\oplus Q_t$
for some projective right module $Q_t$. It implies that
$\Ext^{t+1}(\Tr M,A_A) \simeq \Ext^{1}(\Omega^{t}\Tr M,A_A)\simeq \Ext^{1}(\Tr \mho^t M,A_A)$ and
thus
$\Ext^{1}(\Tr\mho^t M,A_A) \simeq K (\mho^{t} M).$
$\s$
	\medskip

{\bf Remark 1.} {\it For any $t\ge 0$, there is an exact sequence of the form}
$$
  0 @>>>
  \Ext^{t+1}(\Tr M,A_A) @>>> \mho^tM @>\phi_{\mho^tM}>> (\mho^tM)^{**} @>>>
  \Ext^{t+2}(\Tr M,A_A) @>>> 0.
$$
If $t=0$, it is the classical Auslander-Bridger sequence $0 \to \Ext^1(\Tr M,A_A) \to
M @>>> M^{**} \to \Ext^2(\Tr M,A_A)\to 0$ (see [AB], also [ARS]).
	\smallskip
Proof.
The corollary asserts that the kernel of the map $\phi_{\mho^tM}\:\mho^tM \to (\mho^tM)^{**}$
is isomorphic to $\Ext^{t+1}(\Tr M,A_A)$.
On the other hand, the Remark at the end of 2.4 shows that
$\Cok \phi_{\mho^tM} \simeq \K\mho^{t+1}M.$
Since
$\K\mho^{t+1}M \simeq \Ext^{1}(\Tr \mho^{t+1}M,A_A) \simeq
\Ext^{1}(\Omega^{t+1}\Tr M,A_A) \simeq \Ext^{t+2}(\Tr M,A_A)$, it follows that
$\Cok \phi_{\mho^tM} \simeq \Ext^{t+2}(\Tr M,A_A)$.
$\s$
	\medskip
{\bf Remark 2.}  {\it If $M$ is any module, $\mho \Tr M \simeq \Tr \Omega M.$}
	\smallskip
Proof: There is a projective module $P$ such that $\Tr\Tr M \oplus P \simeq M$. According to Lemma
4.4 we have $\mho \Tr M \simeq \Tr\Omega\Tr \Tr M = \Tr\Omega(\Tr\Tr M \oplus P) \simeq \Tr\Omega M.$
 $\s$
	\medskip
{\bf Remark 3.}  In contrast to the isomorphism given in Remark 2, the right modules
$\Omega\Tr M$ and $\Tr\mho M$ discussed in the lemma do not have to be isomorphic.
For example, let $M$ be a module with $M^* = 0$. Then
$\mho M = 0$, thus $\Tr\mho M = 0.$
On the other hand, if $f\:P_1 @>>> P(M)$ is a minimal projective presentation of $M$, then
the kernel of $f^*$ is $M^*$, thus zero, and therefore $\Omega \Tr M \simeq (P(M))^*.$
Thus, we see that the right module $Q$ with $\Omega \Tr M \simeq \Tr \mho M\oplus  Q$
may be non-zero. $\s$
	\medskip
{\bf 4.5. Modules at the end of an $\mho$-path of length $t$.}
	\medskip
{\bf Proposition.} {\it Let $M$ be any module and $t\ge 1.$ The following
conditions are equivalent:
\item{\rm (i)} $\mho^{i-1}M$ is torsionless for $1\le 1\le t$.
\item{\rm (ii)} $M$ is $t$-torsionfree $($thus
$\Ext^i(\Tr M,A_A) = 0$ for $1\le i\le t)$.

If $M$ is indecomposable and not projective, then these conditions are equivalent to
\item{\rm (iii)} $M$ is the end of an $\mho$-path of length $t$.
\par}
	\medskip
Already the special cases $t=1$ and $t=2$ are of interest (but well-known):
{\it A module $M$ is $1$-torsionfree iff $M$ is torsionless} (this is case $t=1$);
{\it a module $M$ is $2$-torsionfree iff both $M$ and $\Omega M$ are torsionless,
thus iff $M$ is reflexive} (this is the case $t=2$, taking into account Corollary 2.4).
These special cases $t=1$ and $t=2$ are discussed at several places; let us refer in
particular to [ARS], Corollary IV.3.3. Our general proof is inspired by [AB].
	\smallskip
Proof of Proposition. For the equivalence of (i) and (ii), see Corollary in 4.4: It asserts
for any $i\ge 1$, that
$\mho^{i-1}M$ is torsionless iff $\Ext^i(\Tr M,A_A) = 0$.

In order to show the equivalence of (i) and (iii), let $M$ be indecomposable and
not projective. If (iii) is satisfied, there is an $\mho$-path of length $t$
ending in $M$. This path has to be $\mho^tM,\ \mho^{t-1}M,\ \dots,\ \mho M,\ M$. This shows
that for any module $\mho^iM$ with $0\le i < t$, there is an arrow starting in
$\mho^iM$, and therefore $\mho^iM$ has to be torsionless.

Conversely, assume that (i) is satisfied. We show (iii) by induction on $t$.
For any $t\ge 1$, there is the arrow $\mho M \to M$, since $M$ is indecomposable, non-projective
and torsionless. According to 3.2, the module $\mho M$ is again indecomposable and non-projective.
Thus, if $t\ge 2$, we can use induction in order to obtain a path of length $t-1$ ending in
$\mho M$, since all the modules $\mho^i(\mho M)$ with $0\le i < t-1$ are torsionless.
$\s$
	\medskip
{\bf 4.6. Proof of Theorem 1.5.}

(1) follows from the fact that $\Ext^t(M,A) = \Ext^{t-1}(\Omega M,A)$ for
$t\ge 2.$ For the special case $t=2$, see Corollary 2.4.
(2) is Proposition 4.5. For (1$'$), (2$'$) and (3), see 4.3.
For (4) and (5), we refer to 4.2(b). Note that in an $\mho$-component of the form
$\Bbb A_n$ with $n\ge 3$, as well as in those of the form $-\Bbb N$, all but precisely two
vertices are the isomorphism classes of reflexive modules, whereas any vertex of
an $\mho$-component of the form $\Bbb N$ is the isomorphism class of a reflexive module.
\hfill$\square$
	\medskip
{\bf 4.7. The adjoint functors $\mho$ and $\Omega.$} Here we collect some important properties of
the construction $\mho$. Some details of the proofs are left to the reader, since the assertions
are not needed in the paper.
\smallskip
If $\Cal C' \subseteq \Cal C$ are full subcategories
of $\mod A$, let $\Cal C/\Cal C'$ be the category with the same objects as $\Cal C$ such that
$\Hom_{\Cal C/\Cal C'}(X,Y)$ is the factor group of $\Hom_{\Cal C}(X,Y)$ modulo the subspace
of all maps $X\to Y$ which factor through a direct sum of modules in $\Cal C'$.
	\smallskip
(1) {\it The functor $\mho$ is  the left adjoint
of the endo-functor $\Omega$ of $\mod A/\add A$.}
Direct verification is easy. But we should also add that
Auslander and Reiten have shown in [AR2, Corollary 3.4]
that the functor $\Tr\Omega\Tr$ is left adjoint to $\Omega$,
and we have identified in 4.4 the functors
$\mho$ and $\Tr\Omega\Tr$.
	\smallskip
(2) Let $\Cal L(A)$ be the full subcategory of all torsionless modules, and $\Cal Z(A)$
the full subcategory of all modules $Z$ with $\Ext^1(Z,A) = 0.$
{\it For any module $M$, the module $\Omega M$ belongs to $\Cal L(A)$,
and the module $\mho M$ belongs to $\Cal Z(A);$ in addition, $\mho M$
has no non-zero projective direct summand.}
	\smallskip
(3) {\it If $Z$
satisfies $\Ext^1(Z, A) = 0$ and has no non-zero projective direct
summand, then $\mho\Omega Z \simeq Z$} (see 3.2).
{\it If $X$ is torsionless and has no non-zero projective direct
summand, then $\Omega\mho X \simeq X$} (see 1.5 or also 3.2). 
 In this way, one shows that
{\it the functors $\Omega$ and $\mho$  provide inverse
categorical equivalences
$$
{\beginpicture
\setcoordinatesystem units <1.5cm,1cm>
\put{$\Cal L(A)/\add(A)$}
 at -.5 0
\put{$\Cal Z(A)/\add(A)$}
 at 2.5 0
\arr{0.5 0.1}{1.5 0.1}
\arr{1.5 -.1}{0.5 -.1}
\put{$\Omega$} at 1 -.3
\put{$\mho$} at 1 .3
\endpicture}
$$}

(4) Thus,
{\it $\Omega$ and $\mho$ provide
inverse bijections between isomorphism classes as follows$:$
$$
{\beginpicture
\setcoordinatesystem units <2cm,1cm>
\put{$\left\{ \matrix \text{\rm indecomposable} \cr
                          \text{\rm non-projective modules $X$}\cr
                   \text{\rm which are torsionless}
                   \endmatrix\right\}$} at -1 0
\put{$\left\{ \matrix \text{\rm indecomposable} \cr
                          \text{\rm non-projective modules $Z$}\cr
                 \text{\rm with  $\Ext^1(Z,A)=0$}
                   \endmatrix\right\}$} at 3 0
\arr{0.5 0.1}{1.5 0.1}
\arr{1.5 -.1}{0.5 -.1}
\put{$\Omega$} at 1 -.3
\put{$\mho$} at 1 .3
\endpicture}
$$
The arrows of the $\mho$-quiver visualize this bijection.}

\medskip
{\bf 4.8. Gorenstein algebras.} Recall that an artin algebra $A$ is said to be
{\it $d$-Gorenstein} provided that the injective dimension of both ${}_AA$ and $A_A$
is equal to $d$. Of course, any algebra of global dimension $d$ is $d$-Gorenstein
The following result of Beligiannis [B2, Proposition 4.4] yields additional examples of
weakly Gorenstein algebras.
	\medskip
{\bf Proposition.} {\it Let $A$ be an artin algebra and assume that
the injective dimension of ${}_AA$ is at most $d$. Then $A$ is right weakly Gorenstein and
any module of the form $\Omega^dM$ is semi-Gorenstein projective.}
\smallskip
Proof: Since the
injective dimension of ${}_AA$ is at most $d$, one knows that
for any module
$M$, the syzygy module $\Omega^d M$ is semi-Gorenstein-projective. [Namely,
for all $i\ge 1$, we have
$\Ext^i(\Omega^d M,A) = \Ext^{d+i}(M,A) = \Ext^i(M,\Sigma^d A) = 0$;
here, $\Sigma N$ denotes the cokernel of an injective envelope of a module $N$.]
This implies that $A$ cannot have any indecomposable module of $\mho$-type $\Bbb N$.
[Namely, if $M$ is of $\mho$-type $\Bbb N$, then
$M$ is $\infty$-torsionfree and therefore $M = \Omega^d(\mho^d M)$. But as
we have seen, this implies that $M$ is semi-Gorenstein-projective, therefore
Gorenstein-projective. Thus $M$ is of $\mho$-type $\Bbb Z$ and not $\Bbb N.$]
Therefore $A$ is right weakly Gorenstein. $\s$
	\medskip
{\bf Corollary 1.}
{\it Let $A$ be $d$-Gorenstein. If an indecomposable non-projective module
$M$ belongs to an $\mho$-path of length $d$, then $M$ is Gorenstein-projective.
If the global dimension of $A$ is  $d$, then there is no
$\mho$-path of length $d$.}
		\smallskip
Proof. Since the $\injdim {}_AA = d$, $A$ is right weakly Gorenstein
and any module $\Omega^dM$ is semi-Gorenstein-projective. Since
$\injdim A_A$ is finite, $A$ is also left weakly Gorenstein, thus the modules $\Omega^d M$ are even
Gorenstein-projective. $\s$
	\medskip
{\bf Corollary 2.} {\it If $A$ is  $d$-Gorenstein, then $A$ has no $\mho$-component
of form $-\Bbb N,\ \Bbb N$ or $\Bbb A_n$ with $n > d$. If the global dimension of $A$ is
$d$, then any $\mho$-component is of form $\Bbb A_n$ with $n\le d.$ }
	\bigskip

{\bf 5. Proof of Theorem 1.4.}
	\smallskip
Since $\add (A) \subseteq   {}^\perp A \subseteq \Cal {F}$, we see that
$\add (A)\subseteq \Cal P(\Cal {F}) = \Cal I(\Cal {F})$.
Thus $\Ext_A^1(X, A) = 0,$ for all $X\in \Cal {F}$.
	\smallskip
For $X\in \Cal {F}$, there is an exact sequence
$0 \longrightarrow K \longrightarrow Q \longrightarrow X \longrightarrow 0$  with $Q\in \Cal {P}(\Cal {F})$ and $K\in \Cal {F}$.
By  $\Cal {P}(\Cal {F}) \subseteq {}^\perp A$ we have $Q\in {}^\perp A$. Thus $\operatorname {Ext}_A^1(X, A) = 0$ and $\operatorname {Ext}_A^{m+1}(X, A) = \operatorname {Ext}_A^m(K, A)$ for $m\ge 1$.
So $\operatorname {Ext}_A^2(X, A) =0,$ and in particular $\operatorname {Ext}_A^2(K, A) =0.$ Repeating this process we see that
$X\in {}^\perp A$. Thus $\Cal {F} \subseteq {}^\perp A$, and hence ${}^\perp A = \Cal {F}$ is Frobenius.
	\smallskip
For $L\in \Cal P(^\perp A)$, consider an exact sequence
$ 0 \longrightarrow K \longrightarrow P \longrightarrow L \longrightarrow 0$  with $P\in \add(A)$. Since $L$ and $P$ are in ${}^\perp A$, $K\in {}^\perp A$. So $\Ext_A^1(L, K) = 0$, thus the exact sequence splits and $L\in \add(A)$. This shows
$\Cal P(^\perp A) \subseteq \add(A) \subseteq \Cal P({}^\perp A)$,
and hence $\Cal P({}^\perp A) = \add(A)$.
	\smallskip
Now consider $X\in  {}^\perp A$.
Since ${}^\perp A$ is Frobenius, there is an exact sequence  $0 \longrightarrow X \longrightarrow I \longrightarrow C \longrightarrow 0$
with $I\in \Cal I(^\perp A) = \Cal P(^\perp A)= \add (A)$ and $C\in {}^\perp A$. So $X$ is torsionless.
This shows that $A$ is left weakly Gorenstein, according to Theorem 1.2.
 $\s$
	\bigskip
{\bf 6. An example.}
\vskip5pt

Let $k$ be a field and $q\in k\setminus\{0\}$.
We consider a $6$-dimensional local algebra
$\Lambda = \Lambda(q)$. If $k$ is infinite, then we show that there are
infinitely many Gorenstein-projective $\Lambda$-modules of dimension 3.
Let $o(q) = |q^{\Bbb Z}|$ be the multiplicative order of $q$. If $o(q)$ is infinite,
we show that there is also
a semi-Gorenstein-projective $\Lambda$-module of dimension 3
which is not Gorenstein-projective.
      \medskip
{\bf 6.1. The algebra $\Lambda = \Lambda(q)$.} The algebra $\Lambda$ is generated by
$x, y, z$, subject to the
relations:
$$
 x^2,\ y^2,\ z^2,\ yz,\ xy+qyx,\ xz-zx,\ zy-zx.
$$
The algebra $\Lambda$ has a basis $1, \ x, \ y, \ z, \ yx,$ and $zx$ and may be visualized
as follows:
$$
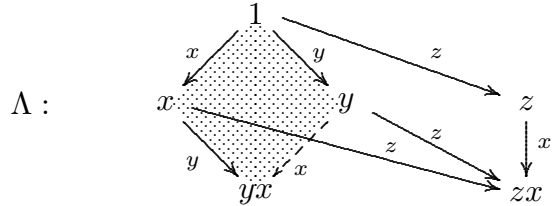
{\beginpicture
      \setcoordinatesystem units <1cm,1cm>
      \put{$\Lambda:$} at -2.5 1
      \put{$1$} at 0 2
      \put{$x$} at -1 1
      \put{$y$} at 1 1
      \put{$yx$} at 0 0
      \arr{0.2 1.8}{0.8 1.2}
      \arr{-.2 1.8}{-.8 1.2}
      \arr{-.8 0.8}{-.2 0.2}
\put{$z$} at 3 1
\put{$zx$} at 3 0
\arr{0.3 1.95}{2.7 1.1}
\arr{1.3 0.9}{2.7 0.15}
\arr{-.7 0.95}{2.7 0.05}
\arr{3 0.8}{3 0.2}
\setdashes <1mm>
\arr{0.8 0.8}{0.2 0.2}
\setlinear
\setshadegrid span <.5mm>
\vshade -1 1 1 <,z,,> 0 0 2 <z,,,> 1 1 1 /
\multiput{$\ssize x$\strut} at -.7  1.6  0.5 0.3 3.2 .6 /
\multiput{$\ssize y$\strut} at 0.7  1.6  -.7 0.4  /
\multiput{$\ssize z$\strut} at 2  1.6  2 0.7  1.5 0.55 /
\endpicture}
$$
Here, we use the following convention: The vertices are the elements of the basis, the arrows are
labeled by $x,y,z$.
A solid arrow $v \rightarrow v'$ labeled say by $x$ means that
$xv = v',$ a dashed arrow $v \dashrightarrow v'$ labeled by $x$
means that $xv$ is a non-zero multiple of
$v'$ (in our case, $xy = -qyx$). If $v$ is a vertex and no arrow starting at $v$ is labeled say by
$x$, then $xv = 0.$

One diamond in the picture has been dotted in order to draw attention to the
relation $xy+qyx$; this relation plays a decisive role when looking at $\Omega M$ for a given
$\Lambda$-module $M$.
      \medskip

We study the following modules $M(\alpha)$ with $\alpha\in k$.
    The module $M(\alpha)$ has a basis $v, v', v''$,
    such that $xv = \alpha v',\ yv = v'$, $zv = v''$, and such that
 $v'$ and $v''$ are annihilated by $x, y, z$. That is,
    $$
    {\beginpicture
    \setcoordinatesystem units <1cm,1cm>
    \put{\beginpicture
    \put{$v$} at 0 1.1
    \put{$v'$} at 0 -.1
    \put{$v''$} at 2 -.1
    \arr{0.08 0.8}{0.08 0.2}
    \arr{.3 .95}{1.7 0}
    \setdashes <1mm>
    \arr{-.08 0.8}{-.08 0.2}
    \put{$\ssize x$\strut} at   -.25 .6
    \put{$\ssize y$\strut} at    .25 .6
    \put{$\ssize z$\strut} at  1 0.7
    \put{$M(\alpha):$} at -1.5 0.8
    \setshadegrid span <.5mm>
    \vshade -.15 0.2 0.81  .15 0.2 0.81 /
    \endpicture} at 0 0
    \endpicture}
    $$
The modules $M(\alpha)$ with $\alpha\in k$
are pairwise non-isomorphic indecomposable $\Lambda$-modules.
    \medskip
For $\alpha\in k,$ we define
$m_\alpha = x - \alpha y\in \Lambda$.
In order to provide a proof of Theorem 1.5, we now collect some general results for the modules
$M(\alpha),\ \Lambda m_\alpha,$ and the right ideals $m_\alpha\Lambda$ which are needed.
	\medskip

{\bf 6.2. The module $M(q)$.}
    \medskip
{\bf Lemma.}  {\it The intersection of the kernels of all the homomorphisms
$M(q) \to {}_\Lambda\Lambda$
is $zM(q) = kv''$ and $M(q)/zM(q) \simeq \Lambda m_1$.
In particular, $M(q)$ is not torsionless and $M(q)^* = (\Lambda m_1)^*$.
}
    \medskip
Proof.  Let $f\:M(q)=\Lambda v \to {}_\Lambda\Lambda$ be a homomorphism.
Let $f(v) = c_1x + c_2y+ c_3z+ c_4yx + c_5zx$ with $c_i\in k$.
By $q f(v') = f(xv) = xf(v) = -c_2qyx + c_3zx$
and $f(v') = f(yv) = yf(v) = c_1yx,$ we get
$c_2 = -c_1$ and $c_3 = 0$.
Thus, $f(v) = c_1(x-y) + c_4yx + c_5zx.$ It follows that $f(v'') = f(zv) = zf(v) = 0$.
This shows that $v''$ is contained in
the kernel of any map $f\:M(q)=\Lambda v \to {}_\Lambda\Lambda$.
On the other hand, the homomorphism
$g\:M(q)=\Lambda v \to \Lambda$ given by $g(v) = x-y = m_1$ has kernel
$kv''$. This completes the proof of the first assertion.

The map $g$ provides a surjective map $p\:M(q) \to \Lambda m_1$ and
$p^*\:M(q)^* \to (\Lambda m_1)^*$ is bijective, thus an isomorphism of right $\Lambda$-modules.
\hfill$\square$
	\medskip
{\bf 6.3. The modules $M(\alpha)$ with $\alpha\in k$.}
We consider now the modules $M(\alpha)$ in general, and relate them to the
left ideals $\Lambda m_\alpha$, and to the right ideals $m_\alpha \Lambda$.
Let us denote by $U_\alpha$ the twosided ideal generated by $m_\alpha$, it is
3-dimensional with basis $m_\alpha,\ yx,\ zx.$
Actually, for any $\alpha\in k,$ the right ideal $m_\alpha\Lambda$ is equal to $U_\alpha$
(but we prefer to write $U_\alpha$ instead of $m_\alpha\Lambda$ when we consider it as a
left module).
For $\alpha\neq 1,$ the left ideal $\Lambda m_\alpha$ is equal to $U_\alpha.$
	\medskip
If $M$ is a module and $m\in M$, we denote by $r(m)\: {}_\Lambda\Lambda \to M$ the
right multiplication by $m$ (defined by $r(m)(\lambda) = \lambda m$.
Similarly, if $N$ is a right $\Lambda$-module and $a\in N$, let $l(a)\:\Lambda_\Lambda
\to N$ be the left multiplication by $a$.

We denote by $u_\alpha\:\Lambda m_\alpha \to \Lambda$
and $u'_\alpha\:m_\alpha\Lambda \to \Lambda$ the canonical embeddings.
	\medskip

(1) {\it The right ideal $m_\alpha\Lambda$ is $3$-dimensional (and equal to $U_\alpha$),
for all $\alpha\in k$.}
	\medskip
(2) {\it The left ideal $\Lambda m_\alpha$ is $3$-dimensional (and equal to $U_\alpha$), for
$\alpha\in k\setminus\{ 1\}$, whereas    $\Lambda m_1$ is 2-dimensional.}
	\medskip
(3) {\it We have $M(\alpha) \simeq \Lambda/U_\alpha$ for all $\alpha\in k$.}
	\smallskip
Proof. The map $r(v):\Lambda \to M(\alpha)$ is surjective (thus a projective cover) and
$$
 r(v)(m_{\alpha}) = m_{\alpha} v = (x-\alpha y)v = xv - \alpha yv = \alpha v' - \alpha v' = 0.
$$
Thus, $\Lambda m_\alpha \subseteq \Ker(r(v))$.
Also, $zx \in \Ker(r(v))$, thus $\Ker(r(v)) = U_\alpha$. This shows that $M(\alpha)$ is isomorphic to $\Lambda/U_\alpha$.
$\s$
	\medskip
(4) {\it For $\alpha\in k\setminus\{1\}$, we have $M(q\alpha) \simeq \Lambda m_\alpha$.}
	\smallskip
Proof.
Consider the map $r(m_\alpha):\Lambda \to \Lambda m_\alpha$.
Since $r(m_\alpha)(m_{q\alpha}) = m_{q\alpha}m_\alpha = 0$,
we see that $U_{q\alpha}\subseteq \Ker(r(m_\alpha))$.
For $\alpha\neq 1,$ the module $\Lambda m_\alpha$ is 3-dimensional, therefore $r(m_\alpha)$
yields an isomorphism $\Lambda/U_{q\alpha} \to \Lambda m_\alpha$.
Using (3) for $M(q\alpha)$, we see that
$
 M(q\alpha) \simeq \Lambda/U_{q\alpha} \simeq \Lambda m_\alpha.
$
$\s$
	\medskip
(5) {\it For any map $f\:\Lambda m_\alpha \to \Lambda,$ there is $\lambda\in \Lambda$
with $f = r(\lambda)u_\alpha$, for all $\alpha\in k$. Thus $u_\alpha$ is a left
$\add(\Lambda)$-approximation.}
	\smallskip
Proof. Let $f\:\Lambda m_\alpha \to \Lambda$ be any map.
Let $f(m_{\alpha}) = c_1x + c_2y+ c_3z+ c_4yx+ c_5zx$ with $c_i\in k$.
Since $f(ym_\alpha) = f(yx)$ and $yf(m_\alpha) = c_1yx$, we see that $f(yx) = c_1yx.$
Since $f(xm_\alpha) = f(-\alpha xy) = q\alpha f(yx) = q\alpha c_1yx$
and $xf(m_\alpha) = c_2xy+c_3zx = -qc_2yx + c_3zx,$ we see that
$q\alpha c_1yx = -qc_2yx + c_3zx$, therefore $c_2 = -\alpha c_1$ and $c_3 = 0.$
Thus,  $f(m_\alpha) = c_1(x-\alpha y)+ c_4yx+ c_5zx$ belongs to
$U_\alpha = m_\alpha\Lambda,$ say $f(m_\alpha) = m_\alpha\lambda$ with $\lambda \in \Lambda$.
Therefore
 $f(m_{\alpha}) = m_\alpha\lambda = r(\lambda)u_\alpha(m_\alpha),$ but this means that
$f = r(\lambda)u_\alpha.$
\hfill$\square$
	\medskip
{\bf 6.4. Lemma.} \ {\it Let $\alpha\in k \setminus\{1\}$. Then there is an $\mho$-sequence}
$$
 0 \to M(q\alpha) \to \Lambda \to M(\alpha)  \to 0.
$$
Proof. According to (3), $M(\alpha) \simeq \Lambda/U_\alpha$. Since $\alpha \neq 1,$
we have $U_\alpha = \Lambda m_\alpha$ by (2). Thus, we have the following exact sequence
$$
 0 \to \Lambda m_\alpha @>u_\alpha>> \Lambda @>>> M(\alpha)  \to 0
$$
 According to (5)
the embedding $u_\alpha\:\Lambda m_\alpha \to \Lambda$ is a left
$\add(\Lambda)$-approximation. Thus, the sequence is an $\mho$-sequence.
Finally, (4) shows that $\Lambda m_\alpha \simeq M(q\alpha).$
 \hfill$\square$
	\medskip
{\bf Corollary 1.} {\it The module $M(0)$ is Gorenstein-projective with $\Omega$-period
equal to $1$.} $\s$

	\medskip
{\bf Corollary 2.} {\it If $o(q) = \infty$, then
the module $M(q)$ is semi-Gorenstein-projective.}
	\smallskip
Proof. We assume that $o(q) = \infty$. Then $q^t \neq 1$ for all $t\ge 1$. By 6.4, all the
sequences
$$
 0 \to M(q^{t+1}) \to \Lambda \to M(q^t)  \to 0.
$$
with $t\ge 1$ are $\mho$-sequences.
They can be concatenated in order to obtain
a minimal projective resolution
of $M(q)$. This shows that $M(q)$ is semi-Gorenstein-projective.
\hfill$\square$

	\medskip
{\bf 6.5. The right $\Lambda$-modules $m_\alpha\Lambda$ and $M(\alpha)^*$.}
We have started in 6.3 to present essential properties of the modules $M(\alpha)$.
We look now also at the modules $m_\alpha\Lambda$ and $M(\alpha)^*$. We
continue the enumeration of the assertions as started in 6.3.
	\medskip
(6) {\it $\Omega (m_{q\alpha}\Lambda) = m_{\alpha}\Lambda$ for all $\alpha\in k$.}
	\smallskip
Proof. We consider the composition of the following right $\Lambda$-module maps
$$
 \Lambda_\Lambda @>l(m_\alpha)>> \Lambda_\Lambda @>l(m_{q\alpha})>> \Lambda_\Lambda
$$
Since $m_{q\alpha}m_{\alpha} = 0,$ the composition is zero.
The image of $l(m_\alpha)$ is the right ideal $m_\alpha\Lambda$, the image of
$l(m_{q\alpha})$ is the right ideal $m_{q\alpha}\Lambda$. Both right ideals are 3-dimensional,
thus the sequence is exact. Thus $m_\alpha\Lambda = \Ker(p)$, for
a surjective map
$p\:\Lambda_\Lambda  \to m_{q\alpha}\Lambda$. Now $p$ is a projective cover, thus
$\Ker (p) = \Omega (m_{q\alpha}\Lambda)$, and therefore
$\Omega(m_{q\alpha}\Lambda) \simeq m_{\alpha}\Lambda$. \hfill$\square$
	\medskip
(7) {\it $(\Lambda m_\alpha)^* = m_\alpha\Lambda$ for all $\alpha\in k.$}
	\smallskip
Proof.
First, let us show that $(\Lambda m_\alpha)^*$ is 3-dimensional. On the one hand,
besides $u_\alpha$, there are
homomorphisms $\Lambda m_\alpha \to \Lambda$ with image $kyx$ and with image $kzx$, which shows that
$(\Lambda m_\alpha)^*$ is at least 3-dimensional.
According to (5), any homomorphism $\Lambda m_\alpha \to \Lambda$ maps into
$\Lambda m_\alpha\Lambda = U_\alpha$.
Since $U_\alpha$ is 3-dimensional, we have $\dim\Hom({}_\Lambda\Lambda,U_\alpha) = 3$, therefore
$\dim (\Lambda m_\alpha)^* = \dim \Hom(\Lambda m_\alpha,\Lambda) = \dim \Hom(\Lambda m_\alpha,U_\alpha)
\le \dim \Hom({}_\Lambda\Lambda,U_\alpha) = 3.$

Second, using again (5), we see that
$(\Lambda m_\alpha)^*$ is, as a right $\Lambda$-module, generated by $u_\alpha$.
Thus, there is a surjective
homomorphism $\theta_\alpha\:\Lambda_\Lambda \to (\Lambda m_\alpha)^*$ of right
$\Lambda$-modules
defined by $\theta_\alpha(1) = u_\alpha$. We have
$$
 (\theta_\alpha(m_{q^{-1}\alpha}))(m_\alpha) = (\theta_\alpha(1)m_{q^{-1}\alpha})(m_\alpha) =
 (u_\alpha m_{q^{-1}\alpha})(m_\alpha) = m_\alpha m_{q^{-1}\alpha} = 0,
$$
therefore $\theta_\alpha(m_{q^{-1}\alpha}) = 0$. It follows that
$\theta_\alpha$ yields a surjective map $\Lambda_\Lambda/m_{q^{-1}\alpha}\Lambda \to (\Lambda m_\alpha)^*.$
Actually, this map has to be an isomorphism, since $m_{q^{-1}\alpha}\Lambda$ is 3-dimensional.
Therefore $\Lambda_\Lambda/m_{q^{-1}\alpha}\Lambda \simeq (\Lambda m_\alpha)^*.$
By (6), $\Lambda_\Lambda/m_{q^{-1}\alpha}\Lambda \simeq m_\alpha\Lambda$.
This completes the proof.
$\s$
	\medskip
(8) {\it $M(q\alpha)^* = m_\alpha\Lambda$ for all $\alpha\in k.$}
	\smallskip
Proof. For $\alpha\neq 1,$ we have $M(q\alpha) \simeq \Lambda m_\alpha$ by (4), thus we use (7).
For $\alpha = 1$, we use 6.2 and then (7).  $\s$
	\medskip
Let us stress that (7) and (8) show that $M(q)^*$ and $(\Lambda m_1)^*$ are isomorphic, namely
isomorphic to $m_1\Lambda$, whereas $M(q)$ and $\Lambda m_1$ themselves are not isomorphic.
	\medskip
(9) {\it Let $\alpha\in k\setminus \{1,q\}.$
For any homomorphism $g\: m_\alpha\Lambda \to \Lambda$ there is $\lambda\in \Lambda$
with $g = l(\lambda)u'_\alpha$. Thus, $u'_\alpha$ is a left
$\add(\Lambda)$-approximation.}
	\smallskip
Proof: Let $g\:m_\alpha\Lambda  \to \Lambda_\Lambda$ be a homomorphism.
We claim that $g(m_\alpha) \in \Lambda m_\alpha$.
Let $g(m_{\alpha}) = c_1x + c_2y+ c_3z+ c_4yx+ c_5zx$ with $c_i\in k$.
Now, $g(m_\alpha x) =  g(-\alpha yx) = -\alpha g(yx)$ and
$g(m_\alpha)x = c_2xy + c_3zx.$
Also,
$g(m_\alpha y) = g(xy) = -qg(yx)$,
and $g(m_\alpha) y = c_1xy+c_3zx = -c_1qyx + c_3zx$,
thus $g(yx) = -q^{-1} g(m_\alpha y) = -q^{-1}( -c_1qyx + c_3zx)
= c_1yx-q^{-1}c_3zx.$
It follows that $c_2yx + c_3zx = -\alpha g(yx) = -\alpha(c_1yx-q^{-1}c_3zx)
  =  -\alpha c_1yx+ \alpha q^{-1}c_3zx.$
Therefore $c_2 = -\alpha c_1$ and $c_3 = \alpha q^{-1}c_3.$ Since we assume that
$\alpha\neq q$, it follows that $c_3 = 0$. Therefore
$g(m_\alpha) = c_1x -\alpha c_1y+ c_3z+ c_4yx+ c_5zx = c_1(x-\alpha y) + c_4yx+ c_5zx $
belongs to $U_\alpha.$
Since we also assume that $\alpha\neq 1$, we have $U_\alpha = \Lambda m_\alpha$.
Thus $g(m_\alpha)\in \Lambda m_\alpha.$

As a consequence, there is $\lambda\in \Lambda$
with $g(m_\alpha) = \lambda m_\alpha$, therefore $g(m_\alpha) = \lambda m_\alpha = l(\lambda) u'_\alpha(m_\alpha).$
It follows that $g = l(\lambda)u'_\alpha.$
$\s$
	\medskip

{\bf 6.6. Lemma.} \ {\it Let $\alpha\in k\setminus \{1,q\}$.
Then there is an $\mho$-sequence of right $\Lambda$-modules}
$$
 0 \to m_{\alpha}\Lambda @>u_\alpha'>> \Lambda_\Lambda @>>> m_{q\alpha}\Lambda  \to 0.
$$
	\smallskip
Proof. This is 6.5 (6) and (9).
$\s$
	\medskip
{\bf 6.7. Proof of Theorem 1.6.}
According to 6.5 (8), we have $M(q)^* = m_{1}\Lambda$.
As we know from 6.2, $M(q)$ is not torsionless.

We assume now that $o(q) = \infty$. The Corollary 2 in 6.4 shows that $M(q)$ is
semi-Gorenstein-projective.
Since $q^{-t} \neq 1$ for all $t\ge 1,$ the sequences
$$
 0 \to m_{q^{-t}}\Lambda @>u_\alpha'>> \Lambda_\Lambda @>>> m_{q^{-t+1}}\Lambda  \to 0
$$
with $t\ge 1$ are $\mho$-sequences, by 6.6. They can be concatenated in order to obtain
a minimal projective resolution
of $m_{1}\Lambda$ and show that $m_{1}\Lambda$ is semi-Gorenstein-projective.
	
Finally, we want to show that $M(q)^{**} = \Omega M(1).$
According to 6.3 (5), the map $u_1\:\Lambda m_1 \to \Lambda$
is a minimal left $\add(\Lambda)$-approximation, thus we may consider as in 2.4 (a)
the following commutative diagram with exact rows:
$$
\CD
  0 @>>> \Lambda m_1    @>u_1>> \Lambda @>\pi_1 >> \Lambda/\Lambda m_1 @>>> 0 \cr
  @.     @VVV                      @|               @VV\phi V \cr
  0 @>>> (\Lambda m_1)^{**} @>>>    \Lambda @>\pi_1^{**}>> (\Lambda/\Lambda m_1)^{**}  @>>>
                                                                    \Ext^1(M'(q)^*,\Lambda_\Lambda) \cr
\endCD
$$
where $\phi = \phi_{\Lambda/\Lambda m_1}.$
The submodule $zx(\Lambda/\Lambda m_1)$ belongs to the kernel of any map $\Lambda/\Lambda m_1) \to \Lambda$,
and it is the kernel of the map $p\:\Lambda/\Lambda m_1 \to M(1)$ defined by $p(\overline 1) = v$.
This shows that $zx(\Lambda/\Lambda m_1)$ is the kernel of $\phi$, thus the image of $\phi$
is just $M(1).$ But the image of $\phi$ coincides with the image of $\pi_1^{**}$.
In this way, we see that $(\Lambda m_1)^{**}$ is the kernel of a projective cover of
$M(1)$, thus equal to $\Omega M(1).$

Of course, $\Omega M(1)$ is decomposable, namely isomorphic to $\Lambda m_1 \oplus kzx$.
$\s$
	\medskip

{\bf 6.8. Proof of Addendum 1.6.} We denote by $q^{\Bbb Z}$ the set of elements of $k$
which are of the form $q^i$ with $i\in \Bbb Z$.
Assume that $\alpha\in k \setminus q^{\Bbb Z}$, thus $q^t\alpha \neq 1$ for all $t\in \Bbb Z$.
According to 6.4, all the
sequences
$$
 0 \to M(q^{t+1}\alpha) \to \Lambda \to M(q^t\alpha)  \to 0
$$
with $t\in \Bbb Z$ are $\mho$-sequences. They can be concatenated in order to obtain
a complete projective resolution
for $M(\alpha)$, thus $M(\alpha)$ is Gorenstein-projective.

The following lemma shows that there are
infinitely many elements $\alpha\in k \setminus q^{\Bbb Z}$.
	\medskip
{\bf Lemma.}
{\it  Assume that $k$ is an infinite field and $q\in k.$ Then $k \setminus q^{\Bbb Z}$ is an
infinite set.}

    \smallskip
We include a proof. The assertion is clear if $o(q)$ is finite. Thus, let
$o(q)$ be infinite (in particular, $q\neq 0$). Assume that
the multiplicative group $k^* = k\setminus \{0\}$ is cyclic, say
$k^* = w^{\Bbb Z}$. Then $o(w)=\infty$, and each element in $k^*$ different from $1$ has infinite multiplicative order. Since $(-1)^2 = 1$, we see that $k$ is of characteristic $2$.
Now $w+1\ne 0$ shows that $w+1 = w^n$ for some $n>1$,
thus $w$ is algebraic over the prime field $\Bbb Z_2$.
Thus $k = \Bbb Z_2(w)$ is a finite field, a contradiction.
Since $k^*$ is not cyclic, there is $a\in k^*\setminus q^{\Bbb Z}$. Then
$a\cdot q^{\Bbb Z}$ is an infinite subset of $k^*\setminus q^{\Bbb Z}$. $\s$
	\bigskip
{\bf 7. Further details for $\Lambda = \Lambda(q)$.}
	\medskip
{\bf 7.1. The $\mho$-components involving modules $M(\alpha)$.}
{\it The only $\mho$-sequences which involve a module of the form $M(\alpha)$ with $\alpha\in k$
are those exhibited in $6.6$.}
\smallskip
Proof. We have to show that there is no $\mho$-sequence ending in $M(1)$ and
no $\mho$-sequence starting in $M(q)$.
Since $\Omega M(1)$ is decomposable, there is no $\mho$-sequence ending in $M(1)$.
By 6.2, the module $M(q)$ is not torsionless, thus no $\mho$-sequence starts in $M(q)$. $\s$
	\medskip
We now want to determine the $\mho$-type of the modules $M(\alpha)$.
According to Corollary 1 in 6.4, $M(0)$ is of $\mho$-type $\widetilde{\Bbb A}_0$. Thus, we now assume
that $\alpha\neq 0.$
	\medskip
{\bf 7.2.} Let us assume that $o(q) = \infty$ (for the case that $o(q)$ is finite, see 7.6).
{\it There are three kinds of $\mho$-components which
involve modules of the form $M(\alpha)$ with $\alpha\in k^*$. There is one component of the form
$-\Bbb N$, it has $M(q)$ as its source, and there is one component of the form
$\Bbb N$, it has $M(1)$ as its sink:
$$
{\beginpicture
\setcoordinatesystem units <1.8cm,1cm>
\put{} at -0.5 0
\put{$M(q^4)$} at 0 0
\put{$M(q^3)$} at 1 0
\put{$M(q^2)$} at 2 0
\put{$M(q)$} at 3 0
\put{$M(1)$} at 4 0
\put{$M(q^{-1})$} at 5 0
\put{$M(q^{-2})$} at 6 0
\setdashes <1mm>
\arr{-.35 0}{-.65 0}
\arr{0.65 0}{0.35 0}
\arr{1.65 0}{1.35 0}
\arr{2.65 0}{2.35 0}
\arr{4.65 0}{4.3 0}
\arr{5.63 0}{5.38 0}
\arr{6.63 0}{6.38 0}
\setshadegrid span <.5mm>
\vshade -.7 -.5 .5  1.35 -.5 .5 /
\vshade  3.65 -.5 .5  6.7 -.5 .5 /
\endpicture}
$$
The remaining ones containing $M(\alpha)$ with $\alpha\neq 0$ and
$\alpha\notin q^{\Bbb Z}$
are of the form $\Bbb Z:$}
$$
{\beginpicture
\setcoordinatesystem units <1.8cm,1cm>
\put{} at -0.5 0
\put{$M(q^4\alpha)$} at 0 0
\put{$M(q^3\alpha)$} at 1 0
\put{$M(q^2\alpha)$} at 2 0
\put{$M(q\alpha)$} at 3 0
\put{$M(\alpha)$} at 4 0
\put{$M(q^{-1}\alpha)$} at 5 0
\put{$M(q^{-2}\alpha)$} at 6 0
\setdashes <1mm>
\arr{-.35 0}{-.65 0}
\arr{0.6 0}{0.4 0}
\arr{1.6 0}{1.4 0}
\arr{2.6 0}{2.4 0}
\arr{3.65 0}{3.35 0}
\arr{4.55 0}{4.3 0}
\arr{5.58 0}{5.42 0}
\arr{6.58 0}{6.42 0}
\setshadegrid span <.5mm>
\vshade -.7 -.5 .5   6.7 -.5 .5 /
\endpicture}
$$
	
\smallskip
\noindent
The positions of the reflexive modules are shaded.
\medskip
According to Theorem 1.5, there are the following observations concerning
the behavior of the modules $M(\alpha)$ with $\alpha\in k$.
	\smallskip
\item{$\bullet$}{\it The module $M(\alpha)$ is Gorenstein-projective iff
   $\alpha\notin q^{\Bbb Z}.$
	\smallskip
\item{$\bullet$}The module $M(\alpha)$ is not Gorenstein-projective,
  but semi-Gorenstein-projective
  iff $\alpha = q^t$ for some $t\ge 1$.
	\smallskip
\item{$\bullet$}The module $M(\alpha)$ is torsionless iff $\alpha\neq q.$
	\smallskip
\item{$\bullet$}The module $M(\alpha)$ is reflexive iff $\alpha\notin\{q,q^2\}$.
	\smallskip
\item{$\bullet$}The module $M(\alpha)$ is not Gorenstein-projective,
  but $\infty$-torsionfree iff $\alpha = q^t$ for some $t\le 0$.\par}
	\medskip

It seems worthwhile to know the canonical maps $\phi_X\:X \to X^{**}$ for the non-reflexive
modules $X = M(q)$ and $X = M(q^2)$. For $M(q)$ we refer to 6.7: there it is shown
that {\it $M(q)^{**} = \Omega M(1)$ and that the image of $\phi_{M(q)}$ is $\Lambda m_1.$}
	\smallskip
It remains to look at $X = M(q^2).$ {\it The module $M(q^2)^{**}$ is
the submodule $\Lambda m_q+\Lambda z$ of $\Lambda$ and $\phi_{M(q^2)}$
is the inclusion map}
$$
  M(q^2) = \Lambda m_q \longrightarrow \Lambda m_q+\Lambda z = M(q^2)^{**}.
$$

Proof. Since $M(q^2)$ is torsionless, the map
$\phi_{M(q^2)}$ is injective.
There is the following commutative diagram with exact rows:
$$
\CD
  0 @>>> M(q^2)    @>u_q>> \Lambda @>\pi_q >> M(q) @>>> 0 \cr
  @.     @VV\phi_{M(q^2)}V                      @|               @VV\phi_{M(q)}V \cr
  0 @>>> M(q^2)^{**} @>u_q^{**}>>    \Lambda @>\pi_q^{**}>> M(q)^{**}  @>>>
                                                                    \Ext^1(M(q^2)^*,\Lambda_\Lambda)
  \ .\cr
\endCD
$$
As we know already, the image of $\phi_{M(q)}$ and therefore of $\pi_q^{**}$, is $\Lambda m_1.$
Thus the kernel of $\pi_q^{**}$ is the submodule $\Lambda m_q+\Lambda z$ of $\Lambda$. Therefore
$M(q^2)^{**} = \Lambda m_q+\Lambda z$ and $\phi_{M(q^2)}$ is the inclusion map $M(q^2) = \Lambda m_q
\longrightarrow \Lambda m_q+\Lambda z = M(q^2)^{**}.$ $\s$

	\medskip
{\bf 7.3. The $\mho$-components involving right $\Lambda$-modules $m_\alpha\Lambda$.}
{\it The $\mho$-sequences which involve a right $\Lambda$-module of the form $m_\alpha\Lambda$
with $\alpha\in k$
are those exhibited in $6.6$ as well as
$$
 0 \to m_q\Lambda @>\bmatrix u_q\strut\cr h\strut \endbmatrix
  >> \Lambda_\Lambda\oplus\Lambda_\Lambda @>>> \mho(m_q\Lambda) \to 0,
$$
and, for $q\neq 1$,
$$
 0 \to m_1\Lambda @>\bmatrix u_1\cr h' \endbmatrix
  >> \Lambda_\Lambda\oplus\Lambda_\Lambda @>>> \mho(m_1\Lambda) \to 0. \quad
$$
Here, $h\:m_q\Lambda \to \Lambda_\Lambda$ is defined by $h(m_q) = z$, whereas
$h'\:m_1\Lambda \to \Lambda_\Lambda$ is defined by $h'(m_1) = zx$.}
	\medskip
Proof. It is easy to check that the map $\bmatrix u_q\cr h \endbmatrix$ and, for $q \neq 1$,
the map
$\bmatrix u_1\strut\cr h'\strut \endbmatrix$ are minimal left $\add(\Lambda_\Lambda)$-approximations.
Clearly, the corresponding cokernels are not torsionless.

In addition, we have to show that there is no $\mho$-sequence ending in $m_{q^2}\Lambda$
or in $m_{q}\Lambda$. But this follows from the fact that the inclusion maps
$u'_q\: m_q\Lambda = \Omega(m_{q^2}\Lambda) \to P(m_{q^2}\Lambda)$ and
$u'_1\: m_1\Lambda = \Omega(m_{q}\Lambda) \to P(m_{q}\Lambda)$
are not $\add(\Lambda_\Lambda)$-approximations. $\s$
	\medskip
{\it Let $o(q) = \infty$} (the case that $o(q) < \infty$ will be discussed in 7.6).
{\it There are five kinds of $\mho$-components
involving right
$\Lambda$-modules of the form $m_\alpha\Lambda$ with $\alpha\in k$, namely a component of the form
$\Bbb N$ with $m_{q^2}\Lambda$ as a sink, a component of the form
$-\Bbb N$ with
$\mho(m_1\Lambda)$ as a source,
and a component of the form $\Bbb A_2$ with sink $m_q\Lambda$ and
source $\mho(m_q\Lambda)$:
$$
{\beginpicture
\setcoordinatesystem units <1.8cm,.8cm>
\put{} at -0.5 0
\put{$m_{q^3}\Lambda$} at 0 0
\put{$m_{q^2}\Lambda$} at 1 0
\put{$m_{q}\Lambda$} at 2 0
\put{$m_{1}\Lambda$} at 3 0
\put{$m_{q^{-1}}\Lambda$} at 4 0
\put{$m_{q^{-2}}\Lambda$} at 5 0
\put{$m_{q^{-3}}\Lambda$} at 6 0
\put{$\mho(m_1\Lambda)$} at 2.1 -1
\put{$\mho(m_q\Lambda)$} at 1.1 -1
\setdashes <1mm>
\arr{-.65 0}{-.35 0}
\arr{0.35 0}{0.65 0}
\arr{3.3 0}{3.65 0}
\arr{4.3 0}{4.65 0}
\arr{5.38 0}{5.63 0}
\arr{6.38 0}{6.63 0}
\arr{1.3 -.8}{1.7 -.3}
\arr{2.3 -.8}{2.7 -.3}
\setshadegrid span <.5mm>
\vshade -.7 -.5 .5  1.35 -.5 .5 /
\vshade  3.65 -.5 .5  6.7 -.5 .5 /

\endpicture}
$$
The $\mho$-components containing right
$\Lambda$-modules $m_\alpha\Lambda$ with
$\alpha\in k \setminus q^{\Bbb Z}$
are of the form $\Bbb Z$:
$$
{\beginpicture
\setcoordinatesystem units <1.8cm,1cm>
\put{} at -0.5 0
\put{$m_{q^3\alpha}\Lambda$} at 0 0
\put{$m_{q^2\alpha}\Lambda$} at 1 0
\put{$m_{q\alpha}\Lambda$} at 2 0
\put{$m_{\alpha}\Lambda$} at 3 0
\put{$m_{q^{-1}\alpha}\Lambda$} at 4 0
\put{$m_{q^{-2}\alpha}\Lambda$} at 5 0
\put{$m_{q^{-3}\alpha}\Lambda$} at 6 0
\setdashes <1mm>
\arr{-.6 0}{-.4 0}
\arr{0.4 0}{0.6 0}
\arr{1.4 0}{1.6 0}
\arr{2.4 0}{2.6 0}
\arr{3.4 0}{3.6 0}
\arr{4.4 0}{4.6 0}
\arr{5.4 0}{5.6 0}
\arr{6.4 0}{6.6 0}
\setshadegrid span <.5mm>
\vshade -.7 -.5 .5  6.7 -.5 .5 /
\endpicture}
$$
In addition, there is the  $\mho$-component consisting of the single right
$\Lambda$-modules $m_0\Lambda$, it is of the form $\widetilde{\Bbb A}_0$.}
	\medskip
For the convenience of the reader, the pictures in 7.1 and 7.2 have been arranged
so that
the $A$-duality is respected. Thus, in 7.1, the arrows are drawn from right to left, in 7.2 from
left to right. Also we recall from 6.3 (8) that the $A$-dual of $M(q\alpha)$ is $m_\alpha\Lambda$,
therefore
the position of $m_\alpha\Lambda$ in the pictures 7.2 is the same as the position of $M(q\alpha)$
in 7.1.
	\medskip
{\bf 7.4.}
We complete the description of the behavior of the modules $M(\alpha)$ started in 7.2.
	\medskip
\item{$\bullet$}{\it The module $M(\alpha)$ is not Gorenstein-projective,
  but $M(\alpha)^*$ is semi-Gorenstein-projective,
  iff $\alpha = q^t$ for some $t\le 1$.}
	\smallskip
\item{$\bullet$}{\it The module $M(\alpha)$ is not Gorenstein-projective,
  but $M(\alpha)^*$ is $\infty$-torsionfree,
  iff $\alpha = q^t$ for some $t\ge 3$.}
\smallskip
Proof.  According to 7.2, the module $M(\alpha)$ is Gorenstein-projective iff
$\alpha\notin q^{\Bbb Z}$. Thus, we can assume that $\alpha = q^t$ for some $t\in \Bbb Z$.
According to 6.3 (8), the module $M(q^t)$ is isomorphic to $m_{q^{t-1}}\Lambda$.
The display of the $\mho$-components shows that $m_{q^{t-1}}\Lambda$ is
semi-Gorenstein-projective iff $t-1 \le 0$, thus iff $t \le 1,$ see Theorem 1.5.
Similarly, we see that $m_{q^{t-1}}\Lambda$ is $\infty$-torsionfree
iff $t-1 \ge 2$, thus iff $t \ge 3.$
$\s$
	\medskip

{\bf 7.5.} We have mentioned in 1.7
that one may use the algebra $\Lambda = \Lambda(q)$ with $o(q)=\infty$
in order to exhibit examples of modules $M$ which satisfy precisely two of the three properties (G1), (G2) and (G3):
	\smallskip
{\it \item{\rm(1)} $M = M(q)$ satisfies {\rm(G1)}, {\rm(G2)}, but not {\rm(G3)}.
\item{\rm(2)} $M = M(q^3)$ satisfies {\rm(G1)}, {\rm(G3)}, but not {\rm(G2)}.
\item{\rm(3)} $M = M(1)$ satisfies {\rm(G2)}, {\rm(G3)}, but not {\rm(G1)}.\par}
\smallskip
Proof: For (1): this is the main assertion of Theorem 1.5. For (2): see 7.2 and 7.3.
For (3): according to 7.2, $M(1)$ is reflexive, but
not Gorenstein-projective. According to 6.3(8), we have
$M(1)^* = m_{q^{-1}}\Lambda$, and  $m_{q^{-1}}\Lambda$ is semi-Gorenstein-projective, see
7.3.
$\s$
	\medskip
Let us look for similar examples for $\Lambda^{\op}$, thus,
for right $\Lambda$-modules $N$.
	\smallskip
{\it \item{\rm(1*)} There is {\bf no} right $\Lambda$-module of the form
   $N = m_\alpha\Lambda$ satisfying {\rm(G1)}, {\rm(G2)}, but not {\rm(G3)}.
\item{\rm(2*)} The right $\Lambda$-module $N = m_{q^{-2}}\Lambda$ satisfies {\rm(G1)},
{\rm(G3)},
   but not {\rm(G2)}.
\item{\rm(3*)} The right $\Lambda$-module  $N = m_{q^2}\Lambda$ satisfies {\rm(G2)},
{\rm(G3)}, but not {\rm(G1)}.
\par}
\smallskip
Proof:
(2*) There starts an infinite $\mho$-path at $N = m_{q^{-2}}\Lambda$, thus
   $N$ satisfies (G1). There ends an $\mho$-path of length 2 at
   $N$, thus  $N$ satisfies (G3).
   Of course, $N^*$ cannot be semi-Gorenstein-projective, since
   otherwise $N$ would be Gorenstein-projective.

(3*) Let $N = m_{q^2}\Lambda$. According to 6.5 (8), $N = M(q^3)^*$. As we know from 7.1,
    $M(q^3)$ is reflexive, thus $N$ is reflexive and $N^* = M(q^3)^{**} = M(q^3)$
    is semi-Gorenstein-projective.

(1*) Assume that $N = m_\alpha\Lambda$ and $N^*$ are both semi-Gorenstein-projective. Since $N$
    cannot be Gorenstein-projective, it is not reflexive. Thus $\alpha \in \{1,q\}$.
    Since $[m_q\Lambda]$ is the sink of an $\mho$-component, $m_\alpha\Lambda$ is not
    semi-Gorenstein-projective. Thus $\alpha = 1.$
    But $(m_1\Lambda)^* = M(q)^{**} = \Omega M(1)$,
    according to  6.5 (8) and Theorem 1.5. As we have mentioned already in the proof 6.7,
    $\Omega M(1) \simeq \Lambda m_1 \oplus k$, where $k$ is the simple $\Lambda$-module.
    We claim that $k$ is not semi-Gorenstein-projective, thus
    $\Omega M(1)$ is not semi-Gorenstein-projective.
	\medskip
{\bf Lemma.} {\it Let $A$ be a local artin algebra which is not self-injective, and $S$ its simple $A$-module. Then $\Ext^i(S,{}_AA) \neq 0$ for all $i\ge 1.$}
\smallskip
Proof: Let $0 \rightarrow \ _AA \rightarrow I_0 \rightarrow I_1 \rightarrow \cdots$ be a minimal injective coresolution. Since $_AA$ is
not injective, all the modules $I_i$ are non-zero. We have $\Ext^i(S, \ _AA)\cong \Hom(S, I_i)$.
$\s$
	\medskip

{\bf 7.6.} Let us look also at the case when $o(q) = n < \infty.$
	\medskip
{\bf Left modules $M(\alpha)$ with $\alpha\in k^*$}
{\it There are two kinds of $\mho$-components which
involve modules of the form $M(\alpha)$ with $\alpha\in k^*$. There is one $\mho$-component of the form
$\Bbb A_n$, it has $M(q)$ as its source, and $M(1)$ as its sink$:$
$$
{\beginpicture
\setcoordinatesystem units <1.85cm,1cm>
\put{$M(1)$} at 0 0
\put{$M(q^{n-1})$} at 1 0
\put{$\cdots$} at 2 0
\put{$M(q^3)$} at 3 0
\put{$M(q^2)$} at 4 0
\put{$M(q)$} at 5 0
\setdashes <1mm>
\arr{0.6 0}{0.3 0}
\arr{1.65 0}{1.4 0}
\arr{2.65 0}{2.35 0}
\arr{3.65 0}{3.35 0}
\arr{4.65 0}{4.35 0}
\setshadegrid span <.5mm>
\vshade -.3 -.5 .5  3.3 -.5 .5 /
\endpicture}
$$
The remaining ones $($containing the modules $M(\alpha)$ with
$\alpha\in k^*\setminus q^{\Bbb Z})$
are directed cycles of cardinality $n:$}
$$
{\beginpicture
\setcoordinatesystem units <1.85cm,1cm>

\put{$M(\alpha)$} at 0 0
\put{$M(q^{n-1}\alpha)$} at 1 0
\put{$\cdots$} at 2 0
\put{$M(q^3\alpha)$} at 3 0
\put{$M(q^2\alpha)$} at 4 0
\put{$M(q\alpha)$} at 5 0
\setdashes <1mm>
\arr{0.6 0}{0.3 0}
\arr{1.65 0}{1.4 0}
\arr{2.65 0}{2.35 0}
\arr{3.65 0}{3.35 0}
\arr{4.65 0}{4.35 0}

\arr{4.7 -.4}{4.8 -.3}
\setquadratic
\plot 0.3 -.38  0.4 -.45  0.6 -.45  2.5 -.45  4.4 -.45  4.6 -.45 4.7 -.38 /
\setlinear
\setshadegrid span <.5mm>
\vshade -.3 -.7 .5  5.3 -.7 .5 /
\endpicture}
$$
\smallskip
\noindent
All modules in the cycles are reflexive. In the $\mho$-component of form $\Bbb A_n$,
the modules $M(q)$ and $M(q^2)$ are not reflexive (they coincide for $o(q) = 1$);
for $o(q)\ge 3$, there are $n-2$ additional
modules $M(1) = M(q^n),\ M(q^{n-2}),\ \cdots, M(q^4),$ $M(q^3)$ in the $\mho$-component,
and these modules are reflexive.
	\medskip
{\bf Right modules $m_\alpha\Lambda$ with $\alpha\in k^*$.}
{\it There is always the $\mho$-component of form $\Bbb A_2$ with
$\Omega(m_q\Lambda)$ as its source and $m_q\Lambda$ as its sink. In addition, for
$n\ge 2$,
there is an $\mho$-component of form $\Bbb A_n$ containing the modules
$m_{q^i}\Lambda$ with $2\le i \le n$ as well as $\Omega(m_1\Lambda);$ it has
$\Omega(m_1\Lambda)$ as its source, and $m_{q^2}\Lambda$ as its sink$:$}
$$
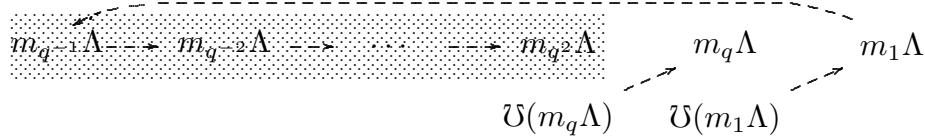
{\beginpicture
\setcoordinatesystem units <1.85cm,.8cm>
\put{} at 5.2 0
\put{$m_{q^{-1}}\Lambda$} at 0 0
\put{$m_{q^{-2}}\Lambda$} at 1 0
\put{$\cdots$} at 2 0
\put{$m_{q^2}\Lambda$} at 3 0
\put{$m_q\Lambda$} at 4 0
\put{$m_{1}\Lambda$} at 5 0
\put{$\mho(m_{1}\Lambda)$} at 4 -1

\put{$\mho(m_{q}\Lambda)$} at 3 -1
\setdashes <1mm>
\arr{0.3 0}{0.6 0}
\arr{1.4 0}{1.65 0}
\arr{2.35 0}{2.65 0}

\arr{3.4 -.6}{3.7 -.3}
\arr{4.4 -.6}{4.7 -.3}
\setshadegrid span <.5mm>
\vshade -.3 -.5 .5  3.3 -.5 .5 /

\arr{0.18 .4}{0.1 .3}
\setquadratic
\plot 0.2 .38  0.3 .5  0.6 .6  2.5 .6  4.4 .6  4.6 .5  4.8 .35 /

\endpicture}
$$
{\it The remaining $\mho$-components $($containing the right modules $m_\alpha\Lambda$ with
$\alpha\in k^*\setminus q^{\Bbb Z})$
are directed cycles of cardinality $n$:}
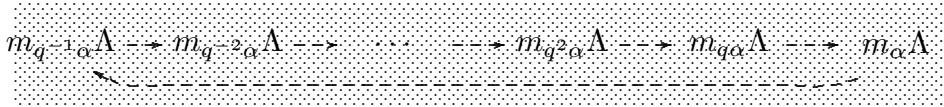
$$
{\beginpicture
\setcoordinatesystem units <1.85cm,1cm>
\put{$m_{q^{-1}\alpha}\Lambda$} at 0 0
\put{$m_{q^{-2}\alpha}\Lambda$} at 1 0
\put{$\cdots$} at 2 0
\put{$m_{q^2\alpha}\Lambda$} at 3 0
\put{$m_{q\alpha}\Lambda$} at 4 0
\put{$m_{\alpha}\Lambda$} at 5 0

\setdashes <1mm>
\arr{0.4 0}{0.6 0}
\arr{1.4 0}{1.65 0}
\arr{2.35 0}{2.65 0}
\arr{3.35 0}{3.65 0}
\arr{4.35 0}{4.65 0}

\arr{0.28 -.38}{0.2 -.3}
\setquadratic
\plot 0.3 -.38  0.4 -.45  0.6 -.45  2.5 -.45  4.4 -.45  4.6 -.45 4.8 -.35 /
\setlinear
\setshadegrid span <.5mm>
\vshade -.3 -.7 .5  5.3 -.7 .5 /
\endpicture}
$$
	\medskip
\noindent
Again, the modules in the cycles are reflexive. In the $\mho$-components of form $\Bbb A_n$
and $\Bbb A_2$,
the modules $m_1\Lambda$ and $\mho(m_1\Lambda)$,
as well as  $m_q\Lambda$ and $\mho(m_q\Lambda)$
are not reflexive; whereas (for $o(q)\ge 3$)
the modules $m_{q^i}\Lambda$ with $2\le i \le n-1$ are reflexive.
	\medskip
Proof: First, let us look at left modules.
According to 7.1, the $\mho$-sequences presented here are the only ones
involving modules of the form $M(\alpha)$. Thus, $[M(q)]$ is a source in the
$\mho$-quiver and $[M(1)]$ is a sink.
This holds true also for $o(q) = 1$: here $q = 1$ and $[M(1)]$ is both a sink and a source, thus a
singleton $\mho$-component (without any arrow).
Finally, for any $n$, the elements $1,q,\dots,q^{n-1}$ are
pairwise different, as are the elements $\alpha, q\alpha,\dots, q^{n-1}\alpha$
for $\alpha\in k \setminus q^{\Bbb Z}$.

For dealing with the right modules, we refer to 7.3.
$\s$
	\medskip
{\bf 7.7.} We have shown in 1.5 that any $\mho$-component is a linearly oriented
quiver of type $\Bbb A_n$ (with $n\ge 1$ vertices),
a directed cycle $\widetilde{\Bbb A}_n$ (with $n+1\ge 1$ vertices),
or of the form $-\Bbb N,$ or $\Bbb N,$ or $\Bbb Z$. Conversely, 7.2 and 7.6 show
that all these cases arise for algebras of the form $\Lambda(q)$.
	\medskip
{\bf 7.8.} A forthcoming paper [RZ] will be devoted to a detailed study of all the
3-dimensional local $\Lambda$-modules for the algebra $\Lambda = \Lambda(q)$.
If $q$ has infinite multiplicative order, we will encounter a whole family of
3-dimensional local modules which are semi-Gorenstein-projective, but not torsionless.
	\bigskip

{\bf 8. Remarks.}
	\smallskip
The first remarks draw the attention to the papers [JS] and [CH]. In 8.1, we show
that the $\Lambda(q)$-modules  $M(q^{-s})$ with $s\ge 0$ and $o(q) = \infty$
satisfy some further conditions which were discussed by Jorgensen and \c Sega.
In 8.2 we show that the algebra $\Lambda(q)$ for $o(q) = \infty$ does not satisfy
the so-called Auslander condition of Christensen and Holm.

In 8.3, we show that essential features of $\Lambda(q)$ are related to
corresponding ones of its subalgebra $\Lambda'(q)$, which is the quantum exterior algebra.
8.4 presents a two-fold covering of $\Lambda(q)$ which has properties
similar to $\Lambda(q)$, but provides for $o(q) = \infty$
examples of semi-Gorenstein-projective modules $M$
which are not Gorenstein-projective, with the additional property that  $\End(M) = k.$
	\medskip
{\bf 8.1. The conditions ($\TR_i$) of Jorgensen and \c Sega.} 
As we have mentioned, Jorgensen and \c Sega have shown in [JS] 
that there exist semi-Gorenstein-projective modules
which are not Gorenstein-projective.
Actually, the main result of [JS] is a stronger assertion.
	\smallskip
Following [JS], we say that an $R$-module $M$ satisfies the condition ($\TR_i$) for some $i\ge 1$
provided $\Ext^i(M,R) = 0$, and that $M$ satisfies
the condition ($\TR_i$) for some $i\le -1$ provided $\Ext^{-i}(\Tr M,R_R) = 0$. Note that ($\TR_i$) is
defined only for $i\neq 0$. Thus, $M$ is semi-Gorenstein-projective if and only if
$M$ satisfies ($\TR_i$) for all $i\ge 1$, and $M$ is $\infty$-torsionfree (i.e., $\Tr M$ is semi-Gorenstein-projective)
if and only if $M$ satisfies ($\TR_i$) for all $i\le -1$. Note that $M$ satisfies $(\TR_i)$ if and only if $\Tr(M)$ satisfies $(\TR_{-i})$. The main theorem of Jorgensen and \c Sega asserts that
{\it there exists a local artinian ring $R$ and a family $M_s$ of $R$-modules, with
$M_s = \Omega M_{s+1}$ for $s\ge 1$, such that $M_s$ satisfies $(\TR_i)$ if and only if $i<s$.}
\smallskip
Such a module $M_s$ satisfies the conditions (G2)
and (G3), and satisfies in addition
the condition that $\Ext^i(M_s,R) = 0$ if and only if $1\le i \le s-1$. Of course, this
is a condition which is much stronger than the negation of (G1).
\smallskip
Let us show that our algebra $\Lambda(q)$ with $o(q)=\infty$ also provides such examples.
Of course, in contrast to the algebra $R$ exhibited by Jorgensen and \c Sega, $\Lambda(q)$ is non-commutative.
There is the following general result:
	
\medskip

{\bf Proposition.} {\it Let $R$ be a local artinian algebra which is not self-injective, with simple $R$-module $S$.

If $M$ is an indecomposable $\infty$-torsionfree module such that $S$ is a proper direct summand of $\Omega M$,
then $M$ satisfies $(\TR_i)$ if and only if $i < 0$.

If $M$ is an indecomposable module such that $M$ satisfies $(\TR_i)$ if and only if $i < 0$,
then for every $s\ge 1$, the module $\mho^{s-1}M$ satisfies $(\TR_i)$ if and only if $i < s$.}
\smallskip
Proof. First, let $M$ be indecomposable, $\infty$-torsionfree, with $\Omega M \cong S\oplus X$ for some non-zero module $X$. Since $M$ is $\infty$-torsionfree,
$M$ satisfies $(\TR_i)$ for $i\le -1$. Since $\Omega M$ is decomposable, we have $\Ext^1(M, R) \ne 0$, i.e., $M$ does not satisfy $(\TR_1)$. By Lemma 7.5, $\Ext^i(S, R) \ne 0$ for all $i\ge 1$. Thus, for $i\ge 2$ we have
$\Ext^i(M, R) \cong \Ext^{i-1}(\Omega M, R) \cong \Ext^{i-1}(S, R) \oplus \Ext^{i-1}(X, R)\ne 0,$ which means that $M$ does not satisfy $(\TR_i)$.

\vskip5pt

Next, assume that $M$ is an indecomposable module such that $M$ satisfies $(\TR_i)$ if and only if $i\le -1$. For $s\ge 1$ consider the module $M_s = \mho^{s-1}M$. For $i\le -1$, $M_s$ satisfies 
satisfies $(\TR_i)$: in fact, by Lemma 4.4, $\Ext^{-i}(\Tr(M_s), R) = \Ext^{-i}(\Tr(\mho^{s-1}M), R) \cong \Ext^{-i}(\Tr(\Tr\Omega^{s-1}\Tr (M)), R)
\cong \Ext^{-i}(\Omega^{s-1}\Tr (M), R) \cong \Ext^{-i+s-1}(\Tr (M), R) = 0.$ 

If $1\le i \le s-1$, then $s-i\ge 1$ and $\Ext^i(M_s, R) = \Ext^i(\mho^{s-1}M, R) \cong \Ext^1(\mho^{s-i}M, R) = 0,$ since $s-i-1\ge 0$ and $\mho^{s-i-1}M$ is torsionless. 

If $i \ge s$, then $i- s +1\ge 1$ shows that  $\Ext^i(M_s,R) \simeq \Ext^{i-s+1}(M, R)\ne 0$, i.e.,
$M_s$ does not satisfy $(\TR_i)$.
$\s$
	\medskip
Application: Let $R = \Lambda = \Lambda(q)$ with $o(q)=\infty.$ Then  $M = M(1)$ is an indecomposable $\infty$-torsionfree module and $S$ is a proper direct summand of $\Omega M$, thus 
the Proposition above shows that for $s\ge 1$, $M_s = \mho^{s-1}M = M(q^{-(s-1)})$ satisfies $(\TR_i)$ if and only if $i<s$.
$\s$
	\medskip
{\bf 8.2. The Auslander condition of Christensen and Holm.} Christensen and Holm [CH]
say that a left-noetherian ring $A$ satisfies the {\it Auslander condition,} provided that
for every finitely generated left $A$-module $M$, there is an integer $b(M)$ with the
following property: if $M'$ is a finitely generated left $A$-module,
then the vanishing $\Ext^{\gg 0}(M,M') = 0$ implies that $\Ext^{> b(M)}(M,M') = 0$.
We are indebted to Christensen and Holm for having drawn our attention to Theorem C of
[CH] which asserts:
{\it If $A$ is a finite-dimensional $k$-algebra $A$ satisfying the Auslander condition, then
A is left weakly Gorenstein} (here, we have taken
into account that a finite-dimensional $k$-algebra has a
dualizing complex, see 3.4 in [CH]). This shows that
{\it the algebra $\Lambda(q)$ with $o(q)= \infty$ does not satisfy the
Auslander condition.}
Actually, this can be seen directly, using the following easy
observation.
	\medskip
{\bf Proposition.} {\it Assume that $A$ is a finite-dimensional $k$-algebra which satisfies
the Auslander condition. Let $N_i$ with $i\in \Bbb Z$ be finite-dimensional right $A$-modules with
$\Omega N_i = N_{i-1}$ for all $i$. If at least one of the modules $N_i$ is semi-Gorenstein-projective, then all the modules $N_i$ are semi-Gorenstein-projective,
thus Gorenstein-projective.}
\smallskip
Proof. Note that $A$ satisfies the Auslander condition if and only if
for every finite-dimensional right $A$-module $N$, there is an integer $c(N)$
such that
for every finite-dimensional right $A$-module $N'$, the vanishing $\Ext^{\gg 0}(N',N) = 0$
implies that $\Ext^{> c(N)}(N',N) = 0$ (here, $ c(N) = b(DN)$, where
$D = \Hom(,-,k)$ denotes the $k$-duality).

We assume that $N_0$ is semi-Gorenstein-projective, whereas $N_1$ is not semi-Gorenstein-projective.
Then we must have $\Ext^1(N_1,A_A)\neq 0$. Since $N_0$ is
semi-Gorenstein-projective,
$\Ext^t(N_0,A_A) = 0$ for all $t\ge 1$ and therefore $\Ext^{t+j}(N_j,A_A) = 0$ for all $t\ge 1$
and $j\ge 0.$ In particular, we have $\Ext^{\gg 0}(N_j,A_A) = 0$ for all $j\ge 0.$
Now we use the Auslander condition with $c = c(A_A)$. Since
$\Ext^{\gg 0}(N_{c+1},A_A) = 0$, we  have $\Ext^{c+1}(N_{c+1},A_A) = 0$. On the other hand,
$\Ext^{c+1}(N_{c+1},A_A) \simeq \Ext^1(N_1,A_A) \neq 0.$ This is a contradiction.
$\s$
	\medskip
For our algebra $\Lambda(q)$ with $o(q)= \infty$, let $N_i = m_{q^i}\Lambda$
with $i\in \Bbb Z$. According to 6.5 (6), we have $\Omega N_i = N_{i-1}$. As we know, the
right module $N_0 = m_1\Lambda = M(q)^*$ is semi-Gorenstein-projective, but
not Gorenstein-projective,  see Theorem 1.6. This shows that
$\Lambda(q)$ with $o(q)= \infty$ does not satisfy the Auslander condition.
	\medskip

{\bf 8.3. The quantum exterior algebra $\Lambda' = \Lambda'(q)$ in two variables}
(see, for example [S]).
Let $\Lambda'$ be the $k$-algebra generated by $x,y$ with the relations
$x^2,\ y^2,\ xy+qyx$. It has a basis $1,\ x,\ y,$
and $yx$. We may use the following picture as an illustration:
$$
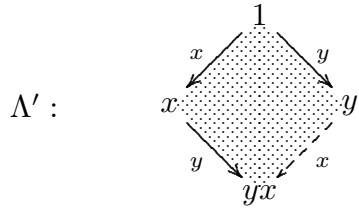
{\beginpicture\setcoordinatesystem units <1cm,1cm>
\put{$1$} at 0 2
\put{$x$} at -1 1
\put{$y$} at 1 1
\put{$yx$} at 0 0
\arr{-.2 1.8}{-.8 1.2}
\arr{0.2 1.8}{0.8 1.2}
\arr{-.8 0.8}{-.2 0.2}
\setdashes <1mm>
\arr{0.8 0.8}{0.2 0.2}
\setshadegrid span <.5mm>
\vshade -1 1 1 <z,z,,> 0 0 2 <z,z,,> 1 1 1 /
\multiput{$\ssize x$\strut} at -.7  1.6  0.7 0.4 /
\multiput{$\ssize y$\strut} at 0.7  1.6  -.7 0.4 /
\put{$\Lambda':$} at -2.5 1
\endpicture}
$$

If we factor out the socle of $\Lambda'$, we obtain the 3-dimensional local algebra
$\Lambda''$ with radical square zero (it is generated by $x,y$ with relations $x^2, y^2, xy, yx$).

Note that $\Lambda'(q)$ is a subalgebra of $\Lambda(q)$, and that
$\Lambda z\Lambda = \Lambda z = \operatorname{span}\{z, zx\}.$
The composition $\Lambda' \hookrightarrow \Lambda \twoheadrightarrow \Lambda/\Lambda z\Lambda$
of the canonical maps is an isomorphism of algebras. In this way, the $\Lambda'$-modules
may be considered as the $\Lambda$-modules which are annihilated by $z$.
We should stress that the elements $m_\alpha = x-\alpha y$ (which play a decisive role
in our investigation) belong to $ \Lambda'$.

For $\alpha\in k$, let $M'(\alpha)$ be the $\Lambda'$-module with basis $v, v'$,
such that $xv =  \alpha v'$,  $yv = v'$, and
$xv' = 0 = yv'$. In addition, we define $M'(\infty)$ as the $\Lambda'$-module
with basis $v, v'$, such that $xv = v', yv = xv' = yv' = 0$.
Here are the corresponding illustrations:
    $$
    {\beginpicture
    \setcoordinatesystem units <1cm,1cm>
    \put{\beginpicture
    \put{$v$} at 0 1.1
    \put{$v'$} at 0 -.1
    \arr{0.08 0.8}{0.08 0.2}
    \setdashes <1mm>
    \arr{-.08 0.8}{-.08 0.2}
    \put{$\ssize x$\strut} at   -.25 .6
    \put{$\ssize y$\strut} at    .25 .6
    \put{$M'(\alpha):$} at -1.5 0.8
    \setshadegrid span <.5mm>
    \vshade -.15 0.2 0.81  .15 0.2 0.81 /
    \endpicture} at 0 0
    \put{\beginpicture
    \put{$v$} at 0 1.1
    \put{$v'$} at 0 -.
    \arr{0 0.8}{0 0.2}
    \put{$\ssize x$\strut} at   -.2 .6
    \put{$M'(\infty):$} at -1.5 0.8
    \endpicture} at 6 0
    \endpicture}
    $$
The modules $M'(\alpha)$ with
$\alpha\in k\cup\{\infty\}$ are pairwise non-isomorphic and indecomposable, and any
two-dimensional indecomposable $\Lambda'$-module is of this form. In particular,
the left ideal $\Lambda' m_\alpha$ is isomorphic to
$M'(q\alpha),$ for any $\alpha\in k\cup\{\infty\}$.
The essential property of the modules $M'(\alpha)$ is the following:
$ \Omega_{\Lambda'} M'(\alpha) = M'(q\alpha).
$
This follows from the fact that $m_{q\alpha}m_\alpha = 0$ and it is this equality which has been
used frequently in sections 6 and 7.
	
For all $\alpha\in k$, $M(\alpha)$ considered as a $\Lambda'$-module,
is equal to $M'(\alpha)\oplus k$,  where $k$ is the simple $\Lambda'$-module.
Also, we should stress that $\rad \Lambda$ considered as a left $\Lambda'$-module is
the direct sum of $I$ and $M'(\infty)$, where $I$ is the indecomposable injective
$\Lambda''$-module.
	\medskip
{\bf 8.4. A variation.} Let $\widetilde\Lambda$ be the algebra defined by a quiver with
two vertices, say labeled by $1$ and $2$, with three arrows $1 \to 2$ labeled by $x,y,z$ and
with three arrows $2 \to 1$, also labeled by $x,y,z$, satisfying the "same" relations
as $\Lambda$ (of course, now we have 14 relations: seven concerning paths $1 \to 2 \to 1$
and seven concerning paths $2 \to 1 \to 2$). Whereas $\Lambda$ is a local
algebra, the algebra $\widetilde\Lambda$ is a connected algebra with two
simple modules $S(1)$ and $S(2)$.

For all the  considerations in
sections 6 and 7, there are corresponding ones for $\widetilde\Lambda$, but always
we have to take into account that now we deal with two simple modules $S(1)$ and
$S(2)$: Corresponding
to the module $M(\alpha)$, there are two different modules,
namely $M^1(\alpha)$ with top
$S(1)$ and $M^2(\alpha)$ with top $S(2)$.
The modules $M^1(\alpha)$ and $M^2(\alpha)$ have similar properties as
$M(\alpha),$ in particular, $M^1(q)$ and $M^2(q)$ are semi-Gorenstein-projective and not
Gorenstein-projective provided that $o(q) = \infty$.
There is one decisive difference between the $\Lambda$-modules and the
$\widetilde\Lambda$-modules: The endomorphism ring of
$M^1(\alpha)$ and $M^2(\alpha)$ is equal to $k$, whereas
the endomorphism ring of any $M(\alpha)$ is 3-dimensional.
	\bigskip
{\bf 9. Questions.}
	\medskip
{\bf 9.1.} We have constructed a module which satisfies the conditions (G1), (G2), but not
(G3). As we have mentioned already in the introduction, it is an open problem whether such
a module does exist in case we deal with commutative rings.
	\medskip

{\bf 9.2.}
One may ask whether or not the finiteness of $\gp A$ implies that $A$ is left weakly  Gorenstein,
There is a weaker question: is $A$ left weakly  Gorenstein, in case all the Gorenstein-projective
$A$-modules are projective?
	\medskip
{\bf 9.3.}
Following Marczinzik [M1, question 1],
one may ask whether a  left weakly Gorenstein artin algebra is also right weakly Gorenstein, thus whether
the existence of an $\mho$-component of the form $\Bbb N$ implies that also an
$\mho$-component of the form $-\Bbb N$ exists.
	\smallskip
Note that if any right weakly Gorenstein algebra is left weakly Gorenstein, then
the Gorenstein symmetry conjecture holds true. Namely, we claim:
{\it If $\injdim {}_AA \le d$ and $\injdim A_A > d$} (the Gorenstein symmetry conjecture asserts
that this should not happen), {\it then $A$ is right weakly Gorenstein, but not left weakly
Gorenstein.}
	\smallskip
Proof. Let $Q$ be an injective cogenerator of $\mod A$.
We assume that $\injdim {}_AA$ is at most $d$. As we have seen in 4.9,
$A$ is right weakly Gorenstein and
any module of the form $\Omega^dM$ is semi-Gorenstein projective.
Now assume that $A$ is also left weakly Gorenstein. Then
all the modules $\Omega^d M$ are Gorenstein-projective. In particular,
$Q' = \Omega^d Q$ is Gorenstein-projective. A well-known argument shows
that if $Q'$ is Gorenstein-projective, then $Q'$ is even projective.
[Namely, assume that $Q'$ is Gorenstein-projective. Then
there is a Gorenstein-projective module $Q''$ such that
$Q' = P'\oplus \Omega^{d+1} Q''$ with $P'$ projective. Now
$
\Ext^1(\Omega^d Q'',Q') \simeq \Ext^{d+1}(Q'',Q') \simeq \Ext^1(Q'',Q) = 0,
$
here the first isomorphism is the usual index shift,
whereas the second comes from the fact that
$Q''$ is (semi-)Gorenstein-projective
and $Q' = \Omega^d Q$
(for a semi-Gorenstein-projective module $N$, and any module $Z$, we have
$\Ext^{i+1}(N,\Omega Z) \simeq \Ext^i(N,Z)$ for all $i \ge 1$).
But $\Ext^1(\Omega^d Q'', P'\oplus \Omega^{d+1} Q'') = 0$ implies that
$\Ext^1(\Omega^d Q'', \Omega^{d+1} Q'') = 0$, thus the canonical exact sequence
$0 \to \Omega^{d+1} Q'' \to P(\Omega^{d} Q'') \to \Omega^d Q''\to 0$ splits and
$\Omega^{d+1} Q''$ has to be projective (even zero). It follows that
$Q' = P'\oplus \Omega^{d+1} Q''$ is projective.]
Since $Q'$ is projective,
the projective dimension of $Q$ is at most $d$. Using duality, we see that
$\injdim A_A \le d$.
$\s$

	\medskip
{\bf 9.4.} Assume that there exists a non-reflexive $A$-module $M$
such that both $M$
and $M^*$ are semi-Gorenstein-projective.  Is then the same true for $A^{\op}$?
Even for $A = \Lambda(q)$ with $o(q) = \infty,$ we do not know the answer.
According to 7.5 (1*), a right $A$-module $N$ of the form
$N = m_\alpha\Lambda(q)$ is reflexive, if both $N$ and $N^*$ are semi-Gorenstein-projective.  But,
there could exist some other right $A$-module $N$ satisfying (G1), (G2) and not (G3).
	\medskip
{\bf 9.5. The Nunke condition.} Does there exist a semi-Gorenstein-projective module $M\neq 0$
with $M^* = 0$\ ? Such a module would be an extreme example of a module
satisfying (G1), (G2) and not (G3). Marczinzik has pointed out that this question
concerns the Nunke condition [H] for $M$,
which asserts that $\Ext^i(M,A) \neq 0$ for some $i\ge 0$, see [J].
Colby and Fuller [CF] have conjectured that the
Nunke condition should hold for any module $M$;
they called this the {\it strong Nakayama conjecture}. The strong Nakayama conjecture
obviously implies the generalized Nakayama conjecture which asserts that
{\it for any simple module $S$ there should exist some $i \ge 0$
such that $\Ext^i(S,A) \neq 0$.} It is known that the Nunke condition is satisfied in case
the finitistic dimension conjecture holds true.

Note that if $M$ is indecomposable and semi-Gorenstein-projective,
then $M^*$ may be decomposable, as Theorem 1.5 shows: the $\Lambda(q)^{\op}$ module $M(q)^*$ is
indecomposable and semi-Gorenstein-projective, but $M(q)^{**}$
is decomposable.
	\medskip
{\bf 9.6. The conditions $(\TR_i)$.} Following
Jorgensen and \c Sega [JS], one may
ask whether an $A$-module which satisfies ($\TR_i$)
for all but finitely many values of $i,$ has to be Gorenstein-projective. In general,
this is not the case, since there is the following proposition.
	\medskip
{\bf Proposition.}
{\it If both $M$ and $M^*$ are semi-Gorenstein-projective, then $M$ satisfies the conditions
$(\TR_i)$ for all $i\notin\{-1,-2\}$.}
\smallskip
Proof. Let $M$ be semi-Gorenstein-projective. Then $M$ satisfies
($\TR_i$) for $i\ge 1$. Since $\Ext^1(M,A) = 0$ for $i=1,2$, Lemma 2.5 asserts that
there is a projective module $Y$ such that $M^* \simeq \Omega^2 \Tr M\oplus Y$.
Assume now that also $M^*$ is semi-Gorenstein-projective. Then for $i\ge 1$,
we have $\Ext^{i+2}(\Tr M,A_A) = \Ext^i(\Omega^2 \Tr M,A_A) = \Ext^i(M^*,A_A) = 0$, thus
$M$ satisfies also ($\TR_i$) for $i\le -3$.  $\s$
	\medskip
Thus, our paper shows that there are
(non-commutative) artinian rings with modules $M$ which satisfy ($\TR_i$)
for all $i\notin\{-1,-2\}$
and which are not Gorenstein-projective. For commutative artinian rings (and this was the setting
considered by Jorgensen and \c Sega) the question is open.
	\bigskip
{\bf 10. References.}
	\medskip
\item {[A1]} M. Auslander, Anneaux de Gorenstein et torsion en alg\'ebre commutative, S\`eminaire d'alg\'ebre
commutative dirig\'e par P. Samuel (1966-1967), tome 1, notes by: M. Mangeney, C. Peskine, L.
Szpiro, Secr\'etariat Math\'ematique, Paris, 2--69(1967).

\item {[A2]} M. Auslander, Coherent functors. In: Proc. the Conf. on
Categorical Algebra. La Jolla 1965. Springer, 189--231.

\item {[AB]}   M. Auslander, M. Bridger, Stable module
theory, Mem. Amer. Math. Soc. 94., Amer. Math. Soc., Providence,
R.I., 1969.

\item{[AR1]}  M. Auslander, I. Reiten, Applications of
contravariantly finite subcategories, Adv. Math. 86(1991), 111--152.

\item{[AR2]}  M. Auslander, I. Reiten, Syzygy modules for Noetherian rings, J. Algebra 183
(1996), 167--185.

\item{[ARS]}  M. Auslander, I. Reiten, S. O. Smal\o, Representation Theory
of Artin Algebras, Cambridge Studies in Advanced Math. 36.
Cambridge University Press, 1995.

\item {[AM]}  L. L. Avramov, A. Martsinkovsky, Absolute, relative, and Tate cohomology of modules of
finite Gorenstein dimension, Proc. London Math. Soc. 85(3)(2002),
393--440.

\item{[B1]}  A. Beligiannis, The homological theory of contravariantly finite subcategories: auslander-buchweitz contexts, gorenstein categories and (co-)stabilization,
Comm. Algebra 28(10) (2000), 4547--4596.

\item{[B2]}  A. Beligiannis, Cohen-Macaulay modules,
(co)torsion pairs and virtually Gorenstein algebras, J. Algebra
288(1)(2005), 137--211.

\item{[B3]} A. Beligiannis, On algebras of finite Cohen-Macaulay type, Adv. Math. 226(2)(2011), 1973--2019.

\item{[Br]} M. Bridger, The $\Ext^i_R(M,R)$ and other invariants of $M$, Brandeis University, Mathematics. Ph.D. (1967).

\item{[Buch]}   R.-O. Buchweitz, Maximal Cohen-Macaulay modules
and Tate cohomology over Gorenstein rings,
Unpublished manuscript, Hamburg (1987), 155pp.

\item{[Che]}  X. W. Chen, Algebras with radical square zero are either self-injective or CM-free,
Proc. Amer. Math. Soc. 140(1)(2012), 93--98.

\item{[Chr]}  L. W. Christensen. Gorenstein Dimensions, Lecture Notes in Math.
1747, Springer-Verlag, 2000.

\item{[CH]} L. W. Christensen, H. Holm, Algebras that satisfy Auslander's condition on
  vanishing of cohomology, Math. Z. 265(2010), 21--40.

\item{[CF]} R. R. Colby, K. R. Fuller, A note on the Nakayama conjecture,
  Tsukuba J. Math. 14(1990), 343--352.

\item{[EJ1]}  E. E. Enochs, O. M. G. Jenda, Gorenstein
injective and projective modules, Math. Z. 220(4)(1995), 611--633.

\item{[EJ2]}  E. E. Enochs, O. M. G. Jenda, Relative
homological algebra, De Gruyter Exp. Math. 30. Walter De Gruyter
Co., 2000.

\item{[H]} D. Happel, Homological conjectures in representation theory of
  finite-dimensional algebras, Unpublished. See:
  https://www.math.uni-bielefeld.de/$\sim$sek/dim2/happel2.pdf (retrieved Aug 6, 2018).

\item{[HH]} C. Huang, Z. Y. Huang, Torsionfree dimension of modules and self-injective
 dimension of rings, Osaka J. Math. 49(2012), 21--35.

\item{[J]} J. P. Jans, Duality in Noetherian rings,
 Proc. Amer. Math. Soc. 12(1961), 829--835.

\item{[JS]}  D. A. Jorgensen, L. M. \c{S}ega, Independence of the
total reflexivity conditions for modules,  Algebras and
Representation Theory 9(2)(2006), 217--226.

\item{[K]} B. Keller,  Chain complexes and stable categories, Manuscripta Math. 67(1990),  379--417.

\item{[L]} T. S. Lam, Lectures on Modules and Rings, Springer, 1999.

\item{[M1]}  R. Marczinzik, Gendo-symmetric algebras, dominant dimensions and Gorenstein
    homological algebra, arXiv:1608.04212.

\item{[M2]}  R. Marczinzik, On stable modules that are not Gorenstein projective, arXiv:1709.01132v3.

\item {[R]}  C. M. Ringel, On the representation dimension of artin algebras,
Bull. the Institute of Math., Academia Sinica, Vol. 7(1)(2012), 33--70.

\item {[RX]} C. M. Ringel, B. L. Xiong, Finite dimensional algebras with Gorenstein-projective nodes.
  In preparation.

\item{[RZ]}  C. M. Ringel, P. Zhang, Gorenstein-projective and semi-Gorenstein-projective modules II,  arXiv:1905.04048.

\item{[S]} S. O. Smal\o, Local limitations of the $\Ext$ functor do not exist, Bull. London Math. Soc. 38(2006), 97--98.

\item{[Y]}  Y. Yoshino, A functorial approach to modules of $G$-dimension zero, Illinois J. Math. 49(3)(2005), 345--367.

\item{[ZX]} P. Zhang, B. L. Xiong. Separated monic representations II: Frobenius subcategories
and RSS equivalences, Trans. Amer. Math. Soc. 372(2)(2019), 981--1021.

\bigskip

\baselineskip=1pt
{\rmk
C. M. Ringel\par
Fakult\"at f\"ur Mathematik, Universit\"at Bielefeld \par
POBox 100131, D-33501 Bielefeld, Germany  \par
ringel\@math.uni-bielefeld.de
\medskip

P. Zhang \par
School of Mathematical Sciences, Shanghai Jiao Tong University \par
Shanghai 200240, P. R. China.\par
pzhang\@sjtu.edu.cn}

\bye